\documentclass[nonblindrev]{informs3} 

\OneAndAHalfSpacedXI

\usepackage{natbib}
\usepackage[ruled]{algorithm2e}
\usepackage{booktabs}
\usepackage{mathrsfs}
\usepackage{subfig}
\usepackage{graphicx}
\usepackage{multicol}
\usepackage{soul}
\usepackage{float}
\usepackage{multirow}
\bibpunct[, ]{(}{)}{,}{a}{}{,}%
\def\bibfont{\small}%
\def\bibsep{\smallskipamount}%

 \newcommand{\x}{{\boldsymbol{x}}}

\newcommand{\R}{\mathbb R}

\newcommand{\dd}{{\boldsymbol{d}}}
\newcommand{\uu}{{\boldsymbol{u}}}
\def\Prob{{\mathbb P}}

\usepackage{subeqnarray, enumitem, amsmath, hyperref, color}

\definecolor{darkblue}{rgb}{0.0,0.0,0.7}
\hypersetup{colorlinks,breaklinks,
	linkcolor=darkblue,urlcolor=darkblue,
	anchorcolor=darkblue,citecolor=darkblue}

\TheoremsNumberedThrough     

\EquationsNumberedThrough    

\begin{document}


\RUNAUTHOR{Cheramin, Cheng, Jiang, and Pan}

\RUNTITLE{Data-Driven Robust Optimization Using ASIUs}

\TITLE{Data-Driven Robust Optimization Using Scenario-Induced Uncertainty Sets}

\ARTICLEAUTHORS{
	\AUTHOR{Meysam Cheramin, Jianqiang Cheng} 
	\AFF{Department of Systems and Industrial Engineering, University of Arizona, Tucson, AZ 85721, USA, \EMAIL{meysamcheramin@email.arizona.edu}, \EMAIL{jqcheng@email.arizona.edu}}
	\AUTHOR{Richard Li-Yang Chen} 
	\AFF{Data Science, Flexport, San Francisco, CA 94102, USA, \EMAIL{rchen@flexport.com}}
	\AUTHOR{Ali Pinar}
	\AFF{Sandia National Laboratories, Livermore, CA 94551, USA,\\\EMAIL{{apinar@sandia.gov}}}
	}

\ABSTRACT{Uncertainty sets are at the heart of robust optimization (RO) because they play a key role in determining the RO models’ tractability, robustness, and conservativeness. Different types of uncertainty sets have been proposed that model uncertainty from various perspectives. Among them, polyhedral uncertainty sets are widely used due to their simplicity and flexible structure to model the underlying uncertainty. However, the conventional polyhedral uncertainty sets present certain disadvantages; some are too conservative while others lead to computationally expensive RO models. This paper proposes a systematic approach to develop data-driven polyhedral uncertainty sets that mitigate these drawbacks. The proposed uncertainty sets are polytopes induced by a given set of scenarios, capture correlation information between uncertain parameters, and allow for direct trade-offs between tractability and conservativeness issue of conventional polyhedral uncertainty sets. To develop these uncertainty sets, we use principal component analysis (PCA) to transform the correlated scenarios into their uncorrelated principal components and to shrink the uncertainty space dimensionality. Thus, decision-makers can use the number of the leading principal components as a tool to trade-off tractability, conservativeness, and robustness of RO models.  We quantify the quality of the lower bound of a static RO problem with a scenario-induced uncertainty set by deriving a theoretical bound on the gap between the optimal value of this problem and that of its lower bound. Additionally, we derive probabilistic guarantees for the performance of the proposed scenario-induced uncertainty sets by developing explicit lower bounds on the number of scenarios required to obtain the desired guarantees. Finally, we demonstrate the practical applicability of the proposed uncertainty sets to trade-off tractability, robustness, and conservativeness by examining a range of knapsack and power grid problems.}


\maketitle

\section{Introduction}

Decision-making in real-world problems is challenging due to the uncertainty involved in them. The challenge is even more significant when the uncertainty is high-dimensional. To overcome this challenge, researchers have proposed various optimization techniques that enable decision-makers to include some knowledge of the uncertainty into their decision-making process to optimize the trade-off between risk and reward. RO is one of these techniques and seeks an optimal solution that is feasible for all realizations within an uncertainty set. RO assumes that all realizations of uncertainty is prescribed by given uncertainty set and hedges against the worst-case scenario in the set \citep{ben1998robust, bertsimas2004price}. 

RO has gained increasing popularity over the last two decades because: (i) it considers uncertainties in the absence of explicit knowledge about their probability distributions; and (ii) its models are usually more tractable than other optimization under uncertainty techniques. Indeed, RO is commonly used in various areas, including but not limited to inventory management, energy management, revenue management, network design, and finance \citep{bertsimas2006robust}. For a detailed review of RO, we refer interested readers to \cite{ben2008selected}, \cite{ben2009robust}, \cite{bertsimas2011theory},  \cite{gabrel2014recent}, and \cite{sozuer2016state}, which provide comprehensive surveys of the RO-related studies.

Uncertainty sets are the core of RO models and play a key role in their performance, greatly impacting solution quality and computational tractability. A well-constructed uncertainty set typically should: (i) capture the most significant aspects of the underlying uncertainty; (ii) be computationally tractable; and (iii) balance robustness and conservativeness of the solution \citep{lorca2014adaptive}. In other words, the uncertainty set should be large enough to include any true realization of uncertainty with high confidence and small enough to exclude pathological scenarios. Since the introduction of RO by \cite{soyster1973convex}, several popular uncertainty sets have been proposed and analyzed. Among them, polyhedral uncertainty sets are the most widely used uncertainty sets due to their computational advantages in deriving linear robust counterparts \citep{lappas2016multi}. Moreover, certain polyhedral uncertainty sets can capture key features of uncertainty, such as asymmetry and correlation, due to their flexibility in including uncertainty data by adjusting their hyperplanes \citep{ning2018data}. The box and
budget uncertainty sets are two popular types of polyhedral uncertainty sets. \cite{soyster1973convex} proposed the box uncertainty set in which each uncertain parameter belongs to a range, while \cite{bertsimas2004price} introduced the budget uncertainty set where the number of the uncertain parameters that are allowed to vary from their nominal values is limited to a pre-specified budget.

Data-driven RO has provided an efficient alternative to traditional decision-making under uncertainty techniques. As a combination of robust and data-driven frameworks, data-driven RO injects a given set of historical data or scenarios into the model through different methods such as constructing a data-driven uncertainty set \citep{bertsimas2006robust}. 

In more recent literature, machine learning techniques have been adopted to develop data-driven uncertainty sets. For instance, \cite{ning2018data, ning2018data2} and \cite{dai2020data} proposed hybrid methods to construct data-driven uncertainty sets by combining the robust kernel density estimation and PCA methods. In other examples,  \cite{shang2017data}, \cite{zhao2019data}, \cite{qiu2019robust}, \cite{shen2020large}, and \cite{mohseni2020data} developed data-driven uncertainty sets using the support vector clustering (SVC) method. Despite the growing popularity of these approaches, there are still many practical limitations. For example, SVC suffers from the curse of dimensionality when uncertainty is high-dimensional \citep{scott2015curse}. 

In this study, we develop data-driven polyhedral uncertainty sets using PCA.  The proposed scenario-induced uncertainty sets have computational benefits for static RO and adaptive robust optimization (ARO) problems by leveraging only a small number of the principal components of uncertainty data while maintaining high solution quality. 

We summarize the key contributions of this paper as follows.
\begin{enumerate}
\item We use PCA to propose a systematic approach for developing data-driven polyhedral uncertainty sets that alleviate the disadvantages of conventional polyhedral uncertainty sets. Unlike the box and budget uncertainty sets, they can capture the correlation information of uncertainty.  Moreover, they are less conservative than the box uncertainty set and computationally cheaper than the convex hull of the uncertainty data. Furthermore, the number of the leading principal components in these uncertainty sets can be used as a tool to trade-off tractability, conservativeness, and robustness of RO models.

\item We quantify the quality of the lower bound of a static RO problem with a scenario-induced uncertainty set by deriving a theoretical bound on the gap between the optimal value of this problem and that of its lower bound. This theoretical bound provides a rough approximation for the optimal value of the static RO Problem, which may not be solved efficiently in practice. Moreover, it determines how many principal components are needed to reach a preferred gap, demonstrating a trade-off between computational burden and solution quality.

\item We provide
probabilistic guarantees for the performance of the proposed uncertainty sets by deriving explicit lower bounds on the number of scenarios required to construct the uncertainty sets with desired probabilistic performance, which complements the existing work. 
\end{enumerate}

The remainder of this paper is organized as follows. In Section \ref{polyhedral}, we provide a concise background on conventional polyhedral uncertainty sets. In Section \ref{Section 3}, after introducing the PCA technique, we propose an efficient approach to construct polyhedral scenario-induced uncertainty sets. In section \ref{ROBoundQuality}, we derive a theoretical bound on the gap between the optimal value of a static RO problem with a scenario-induced uncertainty set and that of its lower bound. In Section \ref{Section 4}, we elaborate on deriving lower bounds on the number of scenario samples required to achieve the desired probabilistic performance guarantees for the developed uncertainty sets.  In Section \ref{ROExperiments}, we conduct extensive computational experiments on RO Knapsack and power grid problems with the proposed uncertainty sets and evaluate their performance.  Finally, Section \ref{Conclusion} contains concluding remarks and future work.

\medskip

\noindent \textbf{Notation}\\
In this paper, we denote scalar values by non-bold symbols, e.g., $m_1$, while we represent vectors by bold symbols in the column form (e.g., $\boldsymbol{u}=\left(u_1,\ldots,u_m\right)^{\top}$ and ${\boldsymbol{w}}$). Similarly, we denote a matrix by bold capital symbols (e.g., $\boldsymbol X$) and indicate its size by $r \times c$, where $r$ and $c$ demonstrate the numbers of rows and columns, respectively. Italic subscripts represent indices, e.g., $c_g$, while non-italic subscripts indicate simplified specifications, e.g., $\mathcal U_\text{box}$. Symbol $||\cdot||$  denotes the Euclidean norm and $|\cdot|$ indicates absolute value. We use $[G]$ to represent the set $\left\{1, 2, \ldots, G\right\}$ for any positive integer number $G$. We reserve symbol $[a,b]$ to represent a range whose minimum and maximum values are $a$ and $b$, respectively. The Euler number is indicated by $e$ while $\boldsymbol e_i$ represents a vector with all zero elements, except for the $i^{th}$ element. Symbol $F(\cdot)$ indicates the cumulative distribution function of a variable and  $\prod$ represents the operator for the product of a sequence. The number of uncertain parameters, i.e., the size of random variable vector, is denoted by $m$ and $\boldsymbol{u}=\left(u_1,\ldots,u_m\right)^{\top} \in \R^{m}$ represents the random variable vector. We adopt $N$ to denote the number of available scenarios for $\boldsymbol{u}$. We reserve symbol $\mathcal S$ to represent the set of the $N$ scenarios, where each scenario is denoted by $\boldsymbol s_j \in \R^{m}$, i.e., $\boldsymbol s_j \in \mathcal S, \ \forall j \in [N]$. The number of utilized principal components in the scenario-induced uncertainty sets is indicated by $m_1$. Symbol $\lceil x\rceil$ represents the smallest integer that is not smaller than $x$. Symbol $\text{unif}(0,1)$ stands for uniform distribution over the interval $[0,1]$.

\section{Polyhedral Uncertainty Sets}\label{polyhedral}

Polyhedral uncertainty sets are widely used in RO because they have a flexible structure to model uncertainty. A general polyhedron uncertainty set is defined as the intersection of closed half-spaces that are represented by linear inequalities of uncertain parameters. More specifically, $\mathcal U_\text{poly}=\{\boldsymbol u: \boldsymbol A_i \boldsymbol u \leq b_i, \ \forall i \in [I]\}$ represents the general formulation of the polyhedral uncertainty set, where $\boldsymbol A_i$ and $b_i$ are the coefficients of its $i^{th}$ linear inequality. 

The box and budget uncertainty sets are two special cases of $\mathcal U_\text{poly}$. \cite{soyster1973convex} introduced $\mathcal U_\text{gen} =\{\boldsymbol u: a_i \leq u_i \leq b_i, \ \forall i \in [m]\}$ as a general box uncertainty set. Alternatively, box uncertainty set can also be defined as follows:
\[ \mathcal U_\text{box}=\{\boldsymbol u: u_i = \bar u_i + z_i \hat u_i,\ -1 \leq z_i \leq 1 , \ \forall i \in [m]\}, \]
where $\bar u_i$ represents the nominal value of $u_i$ and $\hat u_i$ denotes the largest possible deviation of $u_i$, i.e., $u_i$ belongs to range $[\bar u_i - \hat u_i, \bar u_i + \hat u_i]$. \cite{bertsimas2004price} introduced a budget uncertainty set defined as follows:
\[ \mathcal U_\text{budget}=  \{ \boldsymbol u : u_i = \bar u_i + z_i \hat u_i,\ -1 \leq z_i \leq 1 , \  \sum_{i=1}^m |z_i| \le \Gamma, \ \forall i \in [m]\}, \] 
where parameter $\Gamma \in [0,m]$ can be used as a tool to trade-off the conservativeness and robustness of RO models with $\mathcal U_\text{budget}$. Indeed, $\Gamma = 0$ yields the nominal  problem, which does not incorporate any uncertainty, while $\Gamma = m$ results in the most conservative problem, in which $u_i$ is allowed to deviate between its maximum and minimum value. Uncertainty set $\mathcal U_\text{box}$ is a special case of $\mathcal U_\text{budget}$ because $\mathcal U_\text{budget}$ is equivalent to $\mathcal U_\text{box}$ if $\Gamma = m$.     

\begin{figure}[!htb]
	\centering
	\subfloat[\scriptsize{Positively correlated scenarios}]{{\includegraphics[width=5cm, height= 3cm]{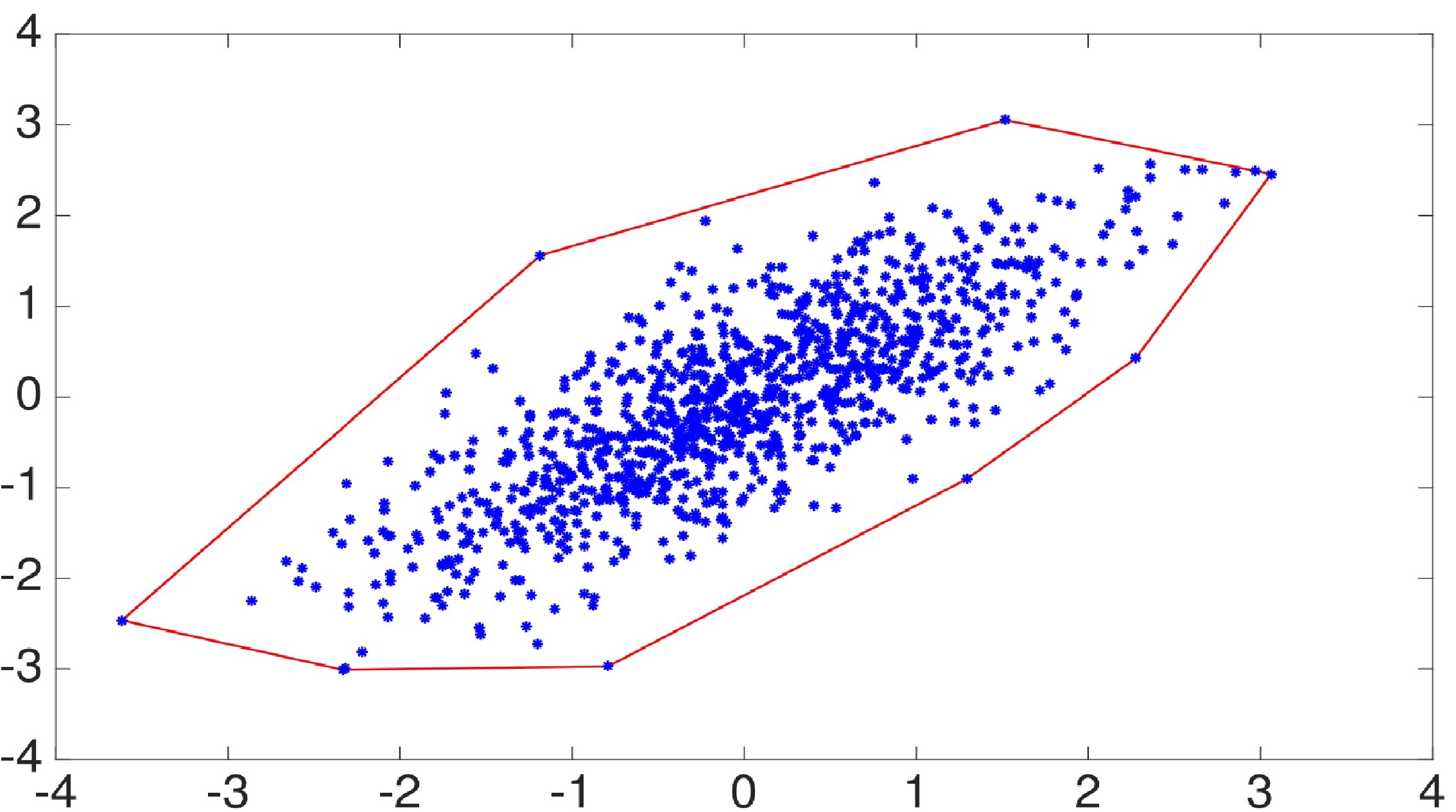} }}%
	\subfloat[\scriptsize{Negatively correlated scenarios}]{{\includegraphics[width=5cm, height= 3cm]{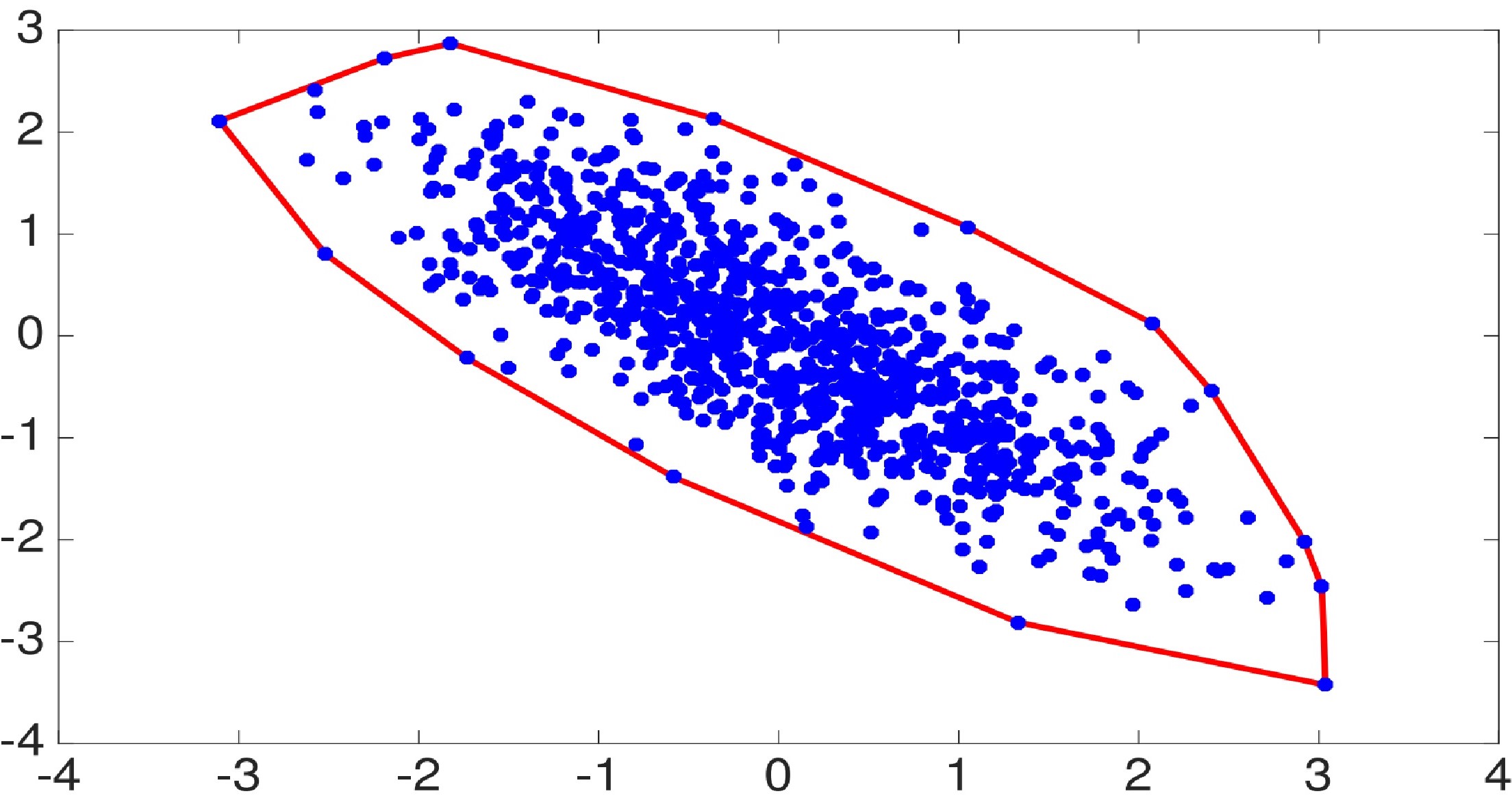} }}%
	\subfloat[\scriptsize{Uncorrelated scenarios}]{{\includegraphics[width=5cm, height= 3cm]{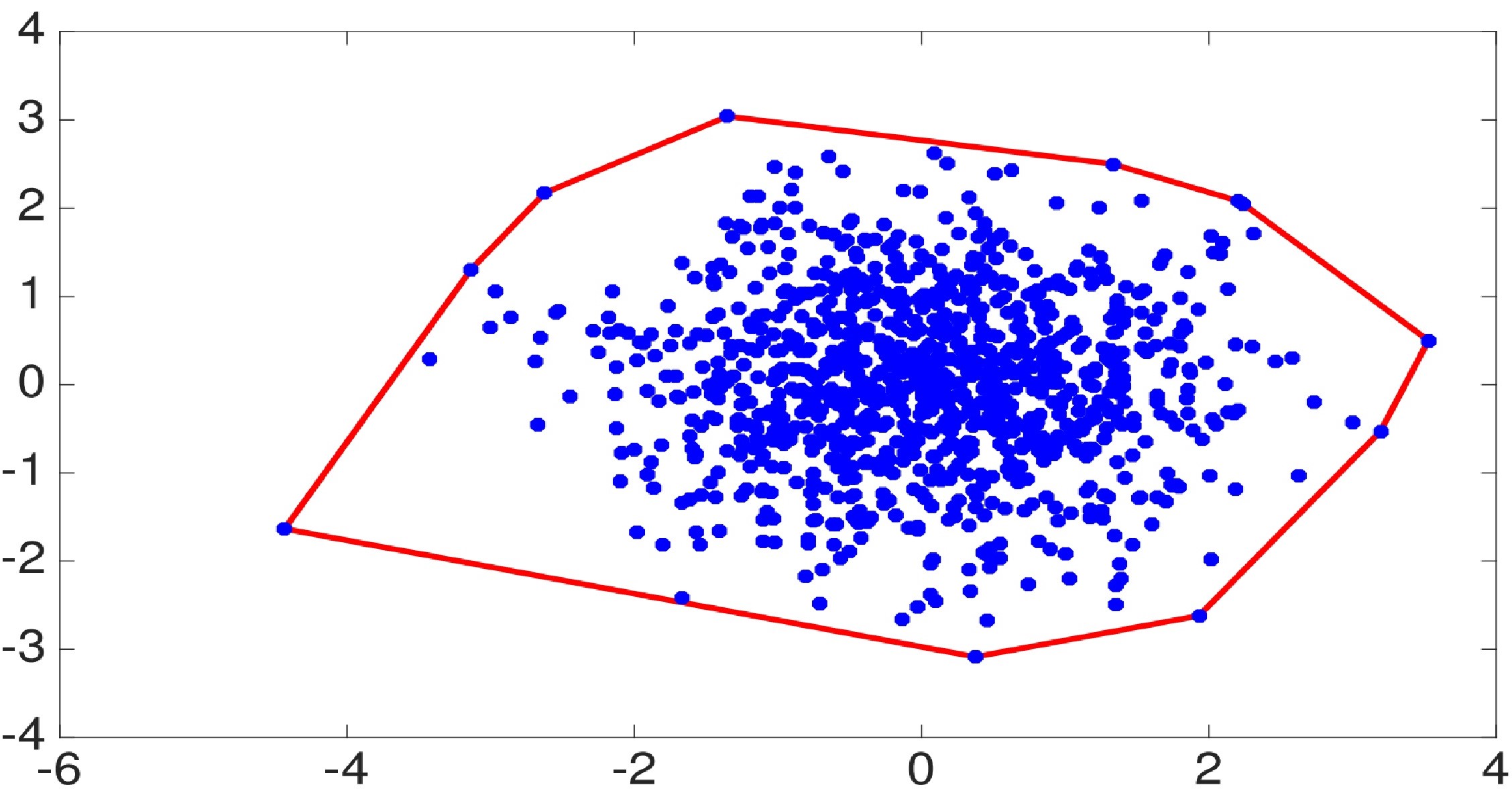} }}%
	\caption{$\mathcal U_{\text{conv}}(\mathcal S)$ of uncertain parameters $u_1, u_2$}
	\label{Fig:siu_sets_orig}%
\end{figure}

\begin{figure}[!htb]
	\centering
	\includegraphics[width=6cm, height= 4cm]{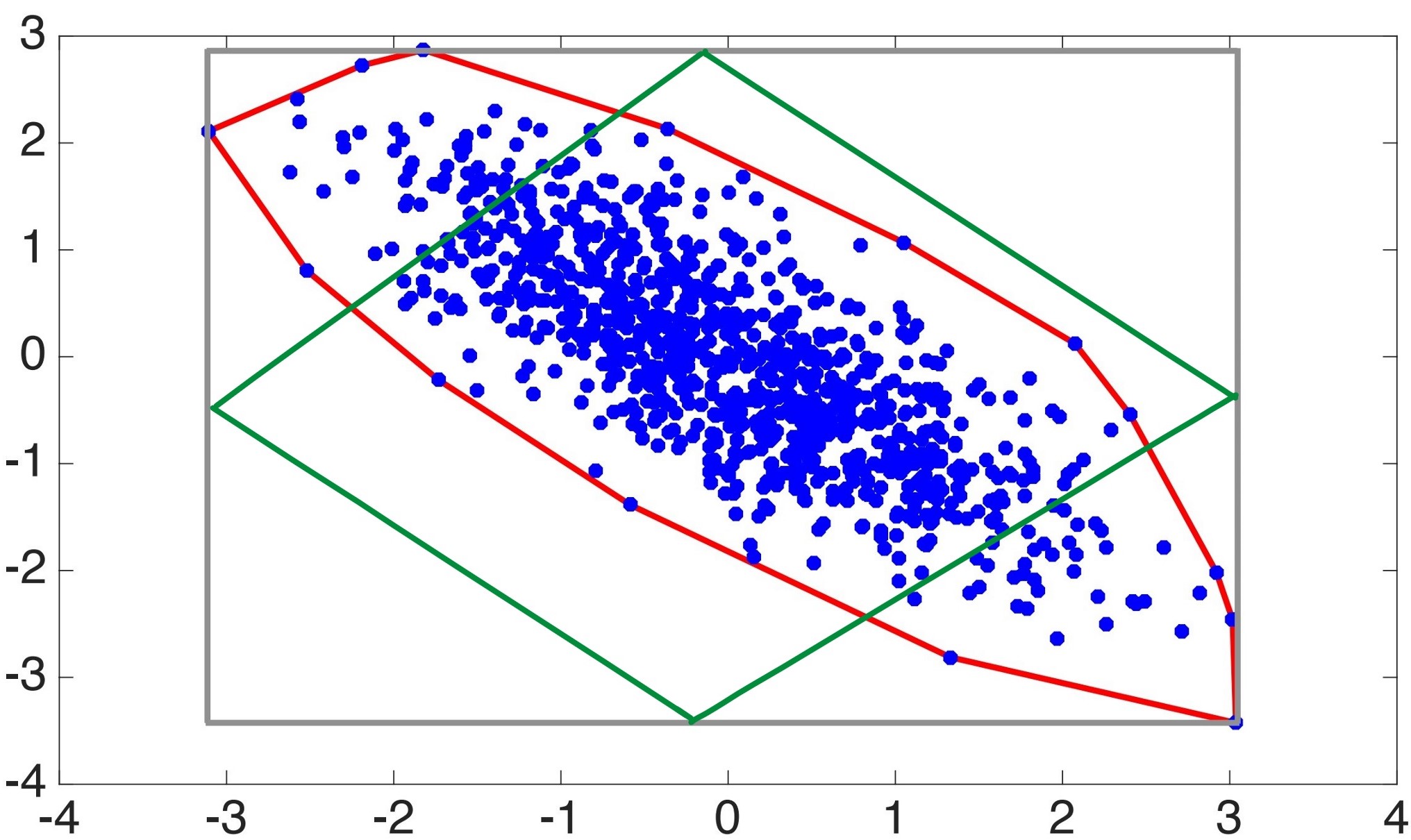}
		\caption{$\mathcal U_\text{box}$, $\mathcal U_\text{budget}$ with $\Gamma=1$, and $\mathcal U_{\text{conv}}(\mathcal S)$ for uncertain parameters $u_1, u_2$ }
	\label{box vs conv(S)}%
\end{figure}

Polyhedral uncertainty sets can be constructed based on historical uncertainty data. These uncertainty sets are referred to as data-driven polyhedral uncertainty sets.  For example, the following convex hull of $\mathcal S$ (i.e., the smallest convex set that includes all $N$ scenarios in $\mathcal S$) can be considered as a scenario-induced polyhedral uncertainty set:
\begin{equation}
	\mathcal U_{\text{conv}}(\mathcal S) = \left \{ \boldsymbol{u}: \  \boldsymbol{u} = \sum_{j = 1}^{N} \alpha_j \boldsymbol s_j, \ \sum_{j = 1}^{N} \alpha_j = 1, \ 0 \leq \alpha_j \le 1, \ \forall j \in [N] \right \}.\nonumber 
\end{equation}
Figure \ref{Fig:siu_sets_orig} illustrates the $\mathcal U_{\text{conv}}(\mathcal S)$ constructed by positively correlated, negatively correlated, and uncorrelated scenarios of $\boldsymbol{u}=(u_1,u_2)^{\top} \in \R^{2}$. Figure \ref{box vs conv(S)} shows $\mathcal U_\text{box}$, $\mathcal U_\text{budget}$ with $\Gamma=1$,  and $\mathcal U_{\text{conv}}(\mathcal S)$ together for this random variable vector. In these figures,  each blue point indicates a  scenario and the gray rectangle, green lozenge, and red polygon represent $\mathcal U_\text{box}$, $\mathcal U_\text{budget}$, and  $\mathcal U_{\text{conv}}(\mathcal S)$, respectively. 

From Figures \ref{Fig:siu_sets_orig} and \ref{box vs conv(S)}, we can observe:  (i) Uncertainty set $\mathcal U_\text{box}$ is the most conservative among them; (ii) Uncertainty sets $\mathcal U_\text{box}$ and $\mathcal U_\text{budget}$ cannot capture the correlation information of uncertainty; and (iii) Uncertainty set $\mathcal U_{\text{conv}}(\mathcal S)$ is the most computationally expensive because it involves many more decision variables  than the other two uncertainty sets ($N$ vs. $m$), leading to larger-size RO formulations. 
Given these observations, it would be highly desirable to develop data-driven polyhedral uncertainty sets that can capture dependent information of uncertainty, alleviate conservatism, and result in more computationally tractable RO models compared to $\mathcal U_{\text{conv}}(\mathcal S)$. To this end, we propose such data-driven polyhedral uncertainty sets in section \ref{Section 3}.

\section{Scenario-Induced Uncertainty Sets}\label{Section 3}

In this section, we develop scenario-induced polyhedral uncertainty sets by leveraging PCA to alleviate the drawbacks of $\mathcal U_\text{box}$, $U_\text{budget}$, and $\mathcal U_{\text{conv}}(\mathcal S)$.  
The merits of the developed uncertainty sets are as follows. First, they explicitly capture the correlation information of uncertainty. 
Second, they yield more tractable RO models compared to $\mathcal U_{\text{conv}}(\mathcal S)$ because their RO formulations are more computationally efficient in comparison with $\mathcal U_{\text{conv}}(\mathcal S)$, due to fewer decision variables. Third, they are less conservative than $\mathcal U_\text{box}$. Fourth, a portion of the principal components of data can be used to improve the tractability and conservativeness of RO models at the expense of robustness reduction.     

In brief, the proposed data-driven approach used to construct scenario-induced uncertainty sets includes the following steps: (i) Calculating the sample mean vector and sample covariance matrix of  $\boldsymbol u$ based on $N$ scenarios; (ii) Obtaining the principal directions of uncertainty data by performing the eigenvalue decomposition on the sample covariance matrix; (iii) Projecting centered scenarios onto each principal direction. In the remainder of this section, we elaborate on these steps in more detail. 

\subsection{Low-rank Approximation with PCA}\label{ROPCA}

The PCA technique enables us to project high-dimensional uncertainty onto a
lower-dimensional space by preserving the components with the highest variance. Moreover, it transforms the correlated uncertain parameters into their uncorrelated principal components \citep{wold1987principal}. We refer interested readers to \cite{wold1987principal} and \cite{reris2015principal} for more information about PCA. 

Let $\bar{\boldsymbol s} =\frac{1}{N} \sum_{j = 1}^{N} \boldsymbol s_j$ be the sample mean of the uncertainty and $\boldsymbol X = [s_{ij}]_{N \times m} $ be the uncertainty data matrix, where the $j^\text{th}$ row represents $\boldsymbol {{s_j}^{\top}} \in \R^{1 \times m}, \ \forall j \in [N]$. Without loss of generality, we center $\boldsymbol X$ at the sample mean by subtracting $\bar{\boldsymbol s}$ from each scenario (row), i.e., $\boldsymbol s_{j0} = \boldsymbol s_j - \bar{\boldsymbol s}, \ \forall j \in [N]$. Therefore, the centered data matrix, denoted by $\boldsymbol X_0$, enables us to approximate the covariance matrix of $\boldsymbol u$ by the sample covariance matrix $\boldsymbol C$ given by $\boldsymbol C = \frac{1}{N-1} \boldsymbol X_0^{\top} \boldsymbol X_0$. 

The PCA technique can be performed by conducting the eigenvalue decomposition (EVD) on $\boldsymbol C$. With $
\boldsymbol X_0^{\top}=\boldsymbol U \boldsymbol \Sigma \boldsymbol V^{\top}$ as the singular value decomposition of $\boldsymbol X_0^{\top}$,  the EAD of $\boldsymbol C$ is as follows:
\begin{equation}\label{decompostion}
	\boldsymbol C = \frac{1}{N-1}\left(\boldsymbol U \boldsymbol \Sigma \boldsymbol V^{\top}\right)\left(\boldsymbol V \boldsymbol \Sigma^{\top} \boldsymbol U^{\top}\right)=\boldsymbol U\left(\frac{\boldsymbol \Sigma \boldsymbol \Sigma^{\top}}{N-1}\right) \boldsymbol U^{\top}=\boldsymbol U \boldsymbol{\Lambda}\boldsymbol U^{\top},\nonumber
\end{equation}
where $\boldsymbol U\in\mathbb{R}^{m\times m}$, $\boldsymbol V\in\mathbb{R}^{N \times N}$,
and $\boldsymbol \Sigma\in\mathbb{R}^{m\times N}$.  The columns of $\boldsymbol U$ and the diagonal entries of $\boldsymbol \Lambda$ represent the eigenvectors and eigenvalues  of $\boldsymbol C$, respectively. The eigenvectors are the principal directions of the centered data, denoted by $\boldsymbol d_i, \ \forall i \in [m]$. The eigenvalue related to each eigenvector represents the variance of the centered data along the corresponding principal direction. Without loss of generality, we assume that the eigenvalues are in non-increasing order. Therefore, the first principal directions characterize most of the variance. The projections of the centered data on the principal directions are called principal components, given by $\boldsymbol s^{'}_{ji}=\frac{\boldsymbol s_{j0} \cdot \boldsymbol d_i }{||\boldsymbol d_i||^2}\boldsymbol d_i, \ \forall j \in [N],  \forall i \in [m]$. 

To reduce the dimensionality of the centered uncertainty data from $m$ to $m_1$, we preserve only the first $m_1$ columns of $\boldsymbol U$ and $m_1 \times m_1$ upper-left entries of $\boldsymbol \Lambda$, which are related to the principal directions with the largest variance. Since the dropped components play the least important role in characterizing the uncertainty, PCA projects the $m$-dimensional uncertainty space onto an $m_1$-dimensional space with the least information loss.

\subsection{PCA-based Polyhedral Uncertainty Sets}\label{PCA-based}

By applying PCA to $\mathcal S$ according to the steps discussed in subsection \ref{ROPCA}, we propose the following scenario-induced uncertainty set:

\begin{footnotesize}
	
	\begin{align}\label{Set 1}
		\mathcal U_{\textsc{pca}}(\mathcal S, m_1) & =   \left \{ \boldsymbol u :   \boldsymbol u = \bar{\boldsymbol s} + \sum_{i=1}^{m_1} \left ( \alpha_i \left(\frac{\overline \omega_i}{||\boldsymbol d_i||}\boldsymbol d_i\right) + (1-\alpha_i)  \left(\frac{\underline \omega_i}{||\boldsymbol d_i||}\boldsymbol d_i\right)\right) + \sum_{i=m_1+1}^m \frac{\overline \omega_i + \underline \omega_i}{2||\boldsymbol d_i||}\boldsymbol d_i, \ 0 \le \alpha_i \le 1, \ \forall i \in [m_1] \right\},\nonumber
	\end{align}
\end{footnotesize}
where 
\begin{small}
	\begin{align}
		\overline \omega_i =\max_{j=1}^{N} \left\{\frac{\boldsymbol s_{j0} \cdot \boldsymbol d_i }{||\boldsymbol d_i||} \right \}\in \R, \quad \text{and} \quad \underline \omega_i =\min_{j=1}^{N} \left\{\frac{\boldsymbol s_{j0} \cdot \boldsymbol d_i}{||\dd_i||} \right \} \in \R, \nonumber
	\end{align}
\end{small}
meaning $(\frac{\overline \omega_i}{||\boldsymbol d_i||}\boldsymbol d_i)$ and $(\frac{\underline \omega_i}{||\boldsymbol d_i||}\boldsymbol d_i)$  are the largest and smallest projected centered scenarios onto the principal direction $\boldsymbol d_i$, respectively. The sample mean $\bar{\boldsymbol s}$ is added to $\mathcal U_{\textsc{pca}}(\mathcal S, m_1)$ because the scenarios have already been centered at $\bar{\boldsymbol s}$.  

\begin{figure}[!htb]
	\centering
	\includegraphics[width=8cm, height= 5cm]{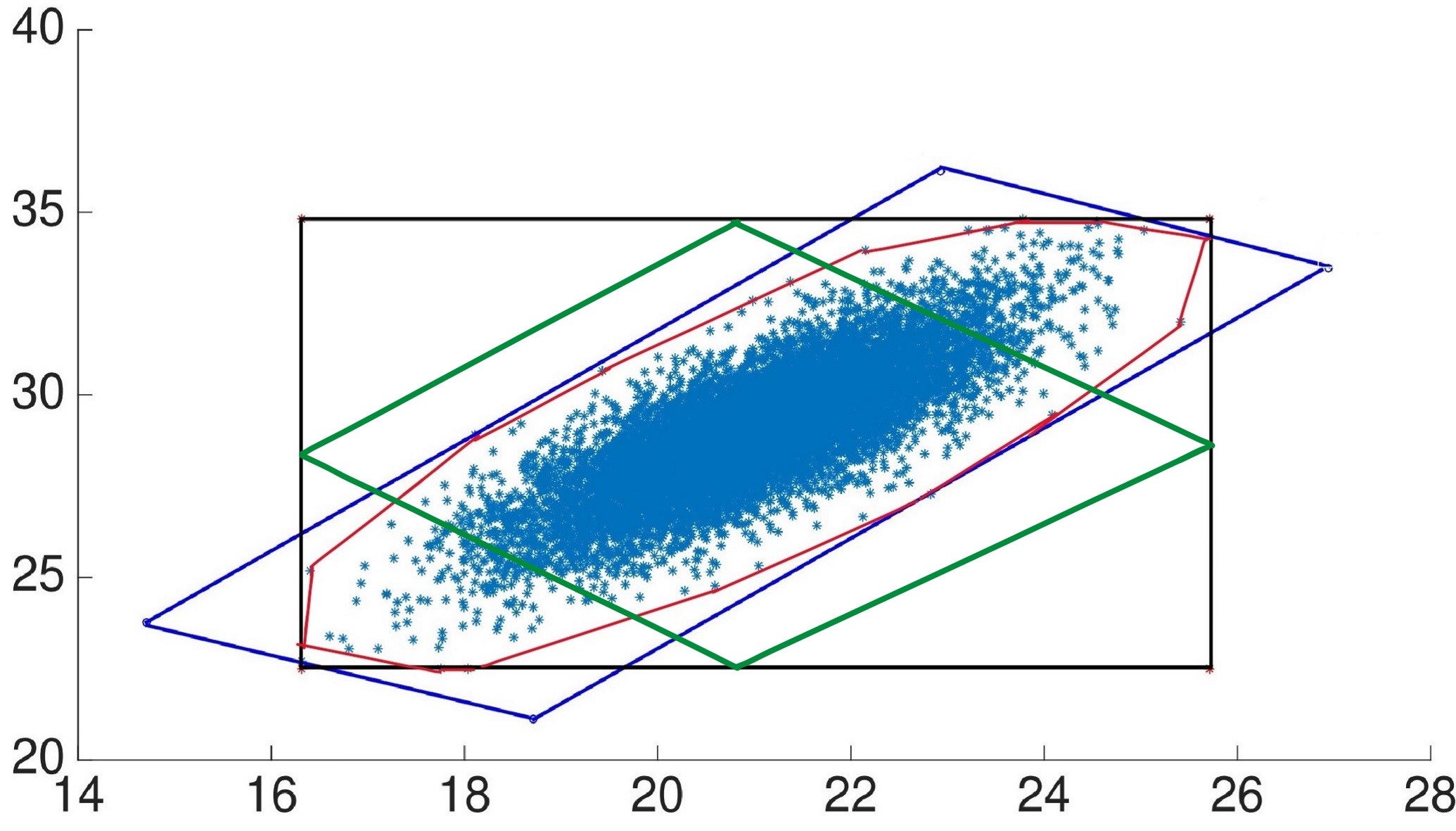}
	\caption{$\mathcal U_{\textsc{pca}}(\mathcal S, 2)$ VS. $\mathcal   U_{\text{conv}}(\mathcal S)$, $\mathcal U_\text{budget}$ with $\Gamma=1$, and $\mathcal U_\text{box}$  for uncertain parameters $u_1, u_2$} \label{box vs conv(S) vs ASIU}%
\end{figure}

In Figure \ref{box vs conv(S) vs ASIU},  the blue rectangle, red polygon, green lozenge, and gray rectangle respectively represent $\mathcal U_{\textsc{pca}}(\mathcal S, 2)$, $\mathcal U_{\text{conv}}(\mathcal S)$, $\mathcal U_\text{budget}$ with $\Gamma=1$, and $\mathcal U_\text{box}$ of $\boldsymbol{u}=(u_1,u_2)^{\top} \in \R^{2}$ for a set of positively correlated scenarios. According to this figure, we have $\mathcal U_{\text{conv}}(\mathcal S) \subseteq  \mathcal U_{\textsc{pca}}(\mathcal S, 2)$, therefore, RO models with $\mathcal U_{\textsc{pca}}(\mathcal S, m_1=m)$ are more conservative (robust) than those with $\mathcal U_{\text{conv}}(\mathcal S)$. On the other hand, $\mathcal U_{\textsc{pca}}(\mathcal S, m_1=m)$ results in more tractable RO models compared to $\mathcal U_{\text{conv}}(\mathcal S)$ because it involves fewer decision variables in RO models in comparison with $\mathcal U_{\text{conv}}(\mathcal S)$.      

\begin{figure}[!htb]
	\centering
	\includegraphics[width=8cm, height= 5cm]{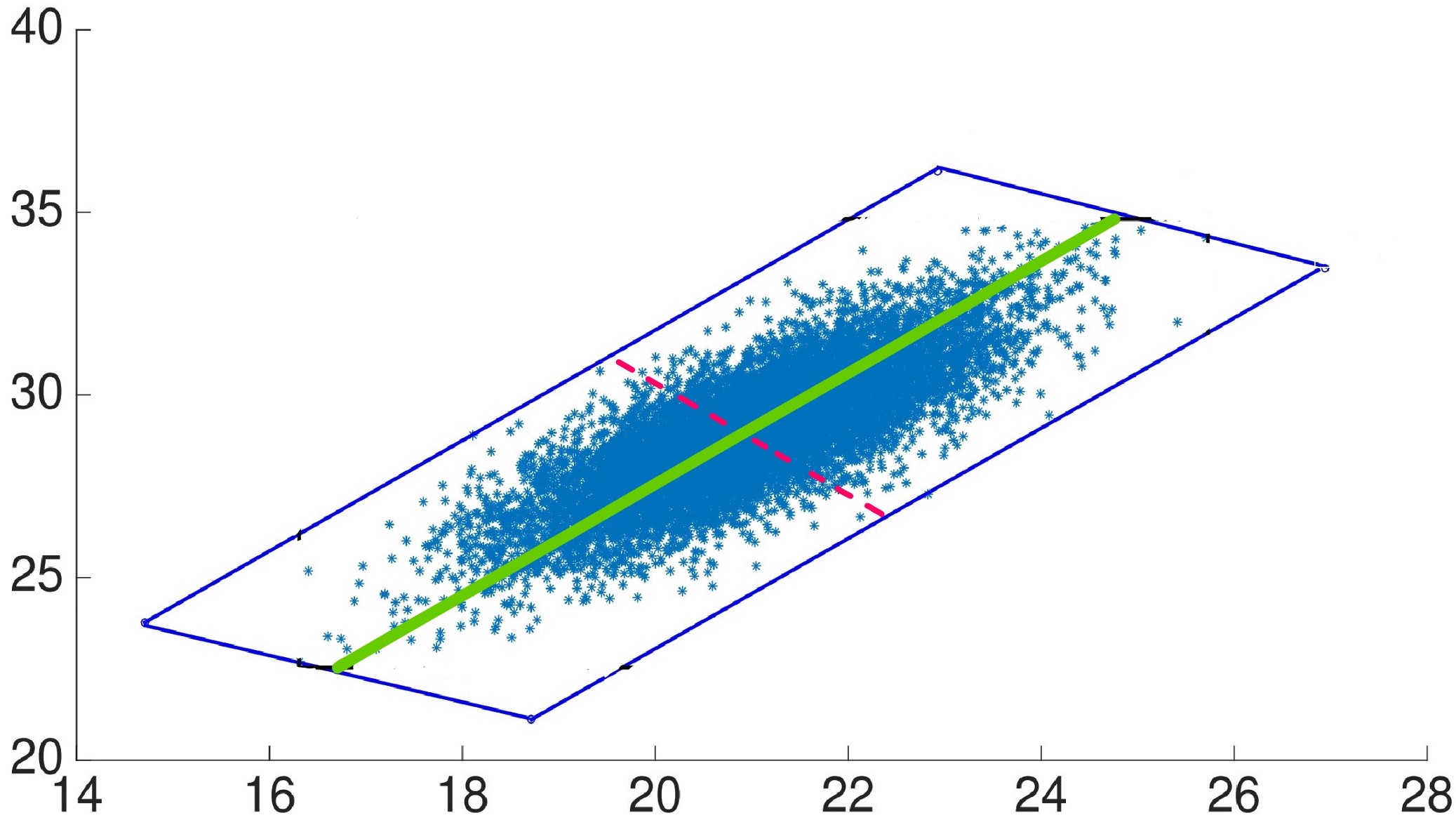}
	\caption{$\mathcal U_{\textsc{pca}}(\mathcal S, 1)$ VS. $\mathcal U_{\textsc{pca}}(\mathcal S, 2)$ for uncertain parameters $u_1, u_2$}
	\label{ASIU(2) VS ASIU(1)}%
\end{figure}

Figure \ref{ASIU(2) VS ASIU(1)} shows $\mathcal U_{\textsc{pca}}(\mathcal S, 2)$ and $\mathcal U_{\textsc{pca}}(\mathcal S, 1)$  for the same set of scenarios. In this example,  $\boldsymbol d_1$ and  $\boldsymbol d_2$, indicated by the green and dashed line respectively, are the principal directions where the most variance of data exists along $\boldsymbol d_1$. The  uncertainty set $\mathcal U_{\textsc{pca}}(\mathcal S, 1)$, which considers $\boldsymbol d_1$ as the only leading principal direction, is the green line segment whose endpoints are generated by $\alpha_1=0$ and $\alpha_1=1$ and any value of $\alpha_1$  between 0 and 1 generates a unique point on this line. As $\mathcal U_{\textsc{pca}}(\mathcal S, 1)$ considers $\boldsymbol d_2$ as the non-leading principal direction, it sets  $\alpha_2=\frac{1}{2}$ to keep only the middle value of $\boldsymbol d_2$ that is located on the green line. Thus, when $m_1$ reduces from 2 to 1, the blue rectangle shrinks to the green line. 

\begin{remark}
	A smaller $m_1$ yields a lower dimensional uncertainty set $\mathcal U_{\textsc{pca}}(\mathcal S, m_1)$, which leads to a more tractable and less conservative (robust) RO model. Therefore, $m_1$ can be used as tool to trade-off tractability and conservativeness of RO models directly.   
\end{remark}

\begin{remark}
	If we choose $\boldsymbol d_i = \boldsymbol e_i, \ \forall i \in [m]$, then the uncertainty set $\hat{\mathcal{U}}_{\textsc{pca}}(\mathcal S, m) =   \{ \boldsymbol u :   \boldsymbol u \in \prod_{i=1}^{m} \left[\bar{\boldsymbol s}+ \underline \omega_i, \ \bar{\boldsymbol s}+ \overline \omega_i \right]\}$ is a box uncertainty set and a special case of $\mathcal U_{\textsc{pca}}(\mathcal S, m)$.
\end{remark}

\begin{remark}
	The intersection of $\mathcal U_{\textsc{pca}}(\mathcal S, m_1)$ and $\hat{\mathcal{U}}_{\textsc{pca}}(\mathcal S, m_1)$, $\mathcal U_\cap(\mathcal S, m_1)=\mathcal U_{\textsc{pca}}(\mathcal S, m_1) \cap \hat{\mathcal{U}}_{\textsc{pca}}(\mathcal S, m_1)$, is another scenario induced uncertainty set.
\end{remark}

\section{Lower Bound Quality  for Static RO}\label{ROBoundQuality}

For any $m_1$ smaller than $m$, a minimization RO problem with $\mathcal U_{\textsc{pca}}(\mathcal S, m_1)$ leads to a lower bound for the same problem with $\mathcal U_{\textsc{pca}}(\mathcal S, m)$. Similarly,  a maximization RO problem with $\mathcal U_{\textsc{pca}}(\mathcal S, m_1)$ results in an upper bound for the same problem with $\mathcal U_{\textsc{pca}}(\mathcal S, m)$. The smaller $m_1$ is, the RO problem with $\mathcal U_{\textsc{pca}}(\mathcal S, m_1)$ yields a looser bound  that is more computationally tractable. In this section, we quantify the quality of the lower bound of a static minimization RO problem with $\mathcal U_{\textsc{pca}}(\mathcal S, m)$ by deriving a theoretical bound
on the gap between its optimal value and the optimal value of its lower bound. Consider the following static RO problem
\begin{equation}\label{RO1}
	Z^*(m) := \min\limits_{\boldsymbol x \in \mathcal X} \max\limits_{\substack{\boldsymbol{u} \in \mathcal U_{\textsc{pca}}\\(\mathcal S, m)}} \ f\left(\boldsymbol{x},\boldsymbol{u}\right), 
\end{equation}
whose lower bound is
\begin{equation}\label{RO6}
	Z^*(m_1) := \min\limits_{\boldsymbol x \in \mathcal X} \max\limits_{\substack{\boldsymbol{u} \in \mathcal U_{\textsc{pca}}\\(\mathcal S, m_1)}} \ f\left(\boldsymbol{x},\boldsymbol{u}\right). 
\end{equation}
The following theorem provides a theoretical bound on the solution quality of problem \eqref{RO6}.
\begin{theorem}\label{ROLBTheoretical0}
	When $f\left(\boldsymbol{x},\boldsymbol{u}\right)$ is piecewise linear convex
	in $\boldsymbol{u}$, i.e., $f\left(\boldsymbol{x},\boldsymbol{u}\right) = \max_{k=1}^K  \left\{y_k^0(\boldsymbol{x})+y_k(\boldsymbol{x})^{\top}\boldsymbol{u}\right\}$
	with both $y_k(\boldsymbol{x})=\left(y_k^1(\boldsymbol{x}),\ldots,y_k^m(\boldsymbol{x})\right)^{\top}$ and $y_k^0(\boldsymbol{x})$ 
	affine in $\boldsymbol{x}$ for any $k \in [K]$, it holds that
	\begin{equation}\label{ROGap1}
		0\leq 	{Z^*}(m)-{Z^*}(m_1)\leq  \max_{k=1}^K  \ \sum_{i=m_1+1}^{m}  |y_k(\hat{\boldsymbol x})^{\top}\overline {\boldsymbol d}_i|  \left (\frac{\overline \omega_i- \underline \omega_i}{2}\right), \nonumber
	\end{equation}
	where $\hat{\boldsymbol x}$ is an optimal solution of the RO problem with $\mathcal U_{\textsc{pca}}(\mathcal S, m_1)$, i.e., Problem \eqref{RO6}, and $\overline {\boldsymbol d_i}=\frac{\boldsymbol d_i}{||\boldsymbol d_i||}, \forall i \in \{m_1+1,\dots,m\}$.
\end{theorem}	

\proof{Proof.} Since Problem \eqref{RO6} is a lower bound of Problem \eqref{RO1}, it is trivial that ${Z^*}_\textup{M}(m)-{Z^*}_\textup{M}(m_1) \geq 0$. In what follows, we derive the upper bound of the gap.  Problem \eqref{RO6} can be rewritten as
\begin{equation}\label{RO2}
	\min\limits_{\boldsymbol x \in \mathcal X} \max\limits_{\substack{\boldsymbol{u} \in  \mathcal U_{\textsc{pca}}\\(\mathcal S, m_1)}} \ \max_{k=1}^K  \ y_k^0(\boldsymbol{x})+y_k(\boldsymbol{x})^{\top}\boldsymbol{u},\nonumber
\end{equation}
which is equivalent to
\begin{equation}\label{RO3}
	\min\limits_{\boldsymbol x \in \mathcal X} \max_{k=1}^K  \ \max\limits_{\substack{\boldsymbol{u} \in \mathcal U_{\textsc{pca}}\\(\mathcal S, m_1)}} y_k^0(\boldsymbol{x})+y_k(\boldsymbol{x})^{\top}\boldsymbol{u}. 
\end{equation}
For clarity, we define $\overline {\boldsymbol d_i}=\frac{\boldsymbol d_i}{||\boldsymbol d_i||}, \ \forall i \in [m]$. With this definition, $\boldsymbol u$ in $ \mathcal U_{\textsc{pca}}(\mathcal S, m_1)$ is rewritten as
\begin{equation}\label{RO4}
	\boldsymbol u = \bar{\boldsymbol s} + \sum_{i=1}^{m_1} \left ( \alpha_i \overline \omega_i\overline {\boldsymbol d_i} + (1-\alpha_i)  \underline \omega_i\overline {\boldsymbol d}_i\right) + \sum_{i=m_1+1}^m \frac{\overline \omega_i + \underline \omega_i}{2}\overline {\boldsymbol d}_i, \ 0 \le \alpha_i \le 1, \ \forall i \in [m_1].
\end{equation}
By plugging \eqref{RO4} to \eqref{RO3}, Problem \eqref{RO3} is reformulated as the following problem:
\begin{equation}\label{RO5}
	\min\limits_{\boldsymbol x \in \mathcal X} \max_{k=1}^K  \ \max\limits_{\substack{0 \le \alpha_i \le 1\\  \forall i \in [m_1]}} y_k^0(\boldsymbol{x})+y_k(\boldsymbol{x})^{\top}\left[\bar{\boldsymbol s} + \sum_{i=1}^{m_1} \left ( \alpha_i \overline \omega_i\overline {\boldsymbol d_i} + (1-\alpha_i)  \underline \omega_i\overline {\boldsymbol d}_i\right) + \sum_{i=m_1+1}^m \frac{\overline \omega_i + \underline \omega_i}{2}\overline {\boldsymbol d}_i\right].\nonumber
\end{equation}
Similarly, Problem \eqref{RO1} is formulated as the following problem:
\begin{equation}\label{RO7}
	\min\limits_{\boldsymbol x \in \mathcal X} \max_{k=1}^K  \ \max\limits_{\substack{0 \le \alpha_i \le 1\\  \forall i \in [m]}} y_k^0(\boldsymbol{x})+y_k(\boldsymbol{x})^{\top}\left[\bar{\boldsymbol s} + \sum_{i=1}^{m} \left ( \alpha_i \overline \omega_i\overline {\boldsymbol d_i} + (1-\alpha_i)  \underline \omega_i\overline {\boldsymbol d}_i\right)\right].\nonumber
\end{equation}
Let $\boldsymbol x^*$ and $\hat{\boldsymbol x}$ be an optimal solution of Problems \eqref{RO1} and \eqref{RO6}, respectively. For clarity, we define
\begin{align}
	h(\boldsymbol x^*,\hat{\alpha}_{m})=&\ y_k^0(\boldsymbol x^*)+y_k(\boldsymbol x^*)^{\top}\left[\bar{\boldsymbol s} + \sum_{i=1}^{m} \left ( \alpha_i \overline \omega_i\overline {\boldsymbol d_i} + (1-\alpha_i)  \underline \omega_i\overline {\boldsymbol d}_i\right)\right], \nonumber \\
	g(\hat{\boldsymbol x},\hat{\alpha}_{m_1})=&\ y_k^0(\hat{\boldsymbol x})+y_k(\hat{\boldsymbol x})^{\top}\left[\bar{\boldsymbol s} + \sum_{i=1}^{m_1} \left ( \alpha_i \overline \omega_i\overline {\boldsymbol d_i} + (1-\alpha_i)  \underline \omega_i\overline {\boldsymbol d}_i\right) + \sum_{i=m_1+1}^m \frac{\overline \omega_i + \underline \omega_i}{2}\overline {\boldsymbol d}_i\right],\nonumber		
\end{align}
where $\hat{\alpha}_{m}=\{\alpha_1, \dots, \alpha_{m}\}$ and $\hat{\alpha}_{m_1}=\{\alpha_1, \dots, \alpha_{m_1}\}$. With the definitions, we have
\[ {Z^*}(m)-	{Z^*}(m_1)=  \max_{k=1}^K  \ \max\limits_{\substack{0 \le \alpha_i \le 1\\  \forall i \in [m]}} h(\boldsymbol x^*,\hat{\alpha}_{m}) - \max_{k=1}^K  \ \max\limits_{\substack{0 \le \alpha_i \le 1\\  \forall i \in [m_1]}} g(\hat{\boldsymbol x},\hat{\alpha}_{m_1}).\]
Since $\hat{\boldsymbol x}$ is a feasible solution of Problem \eqref{RO1} as well, we have
\begin{small}
	\begin{align}
		\Leftrightarrow &	{Z^*}(m)-	{Z^*}(m_1) \leq  	\max_{k=1}^K  \ \max\limits_{\substack{0 \le \alpha_i \le 1\\  \forall i \in [m]}} h(\hat{\boldsymbol x},\hat{\alpha}_{m}) - \max_{k=1}^K  \ \max\limits_{\substack{0 \le \alpha_i \le 1\\  \forall i \in [m_1]}} g(\hat{\boldsymbol x},\hat{\alpha}_{m_1}),\nonumber\\
		\Leftrightarrow &{Z^*}(m)-	{Z^*}(m_1) \leq \max_{k=1}^K  \ \max\limits_{\substack{0 \le \alpha_i \le 1\\  \forall i \in [m]}} y_k^0(\hat{\boldsymbol x})+y_k(\hat{\boldsymbol x})^{\top}\left[\bar{\boldsymbol s} + \sum_{i=1}^{m} \left ( \alpha_i \overline \omega_i\overline {\boldsymbol d_i} + (1-\alpha_i)  \underline \omega_i\overline {\boldsymbol d}_i\right)\right]  \nonumber\\
		& - \max_{k=1}^K  \ \max\limits_{\substack{0 \le \alpha_i \le 1\\  \forall i \in [m_1]}}\left(y_k^0(\hat{\boldsymbol x})+y_k(\hat{\boldsymbol x})^{\top}\left[\bar{\boldsymbol s} + \sum_{i=1}^{m_1} \left ( \alpha_i \overline \omega_i\overline {\boldsymbol d_i} + (1-\alpha_i)  \underline \omega_i\overline {\boldsymbol d}_i\right) + \sum_{i=m_1+1}^m \frac{\overline \omega_i + \underline \omega_i}{2}\overline {\boldsymbol d}_i\right]\right)\nonumber\\
		\Leftrightarrow & {Z^*}(m)-	{Z^*}(m_1) \leq  \max_{k=1}^K  \ \max\limits_{\substack{0 \le \alpha_i \le 1\\  \forall i \in \{m_1+1,\dots,m\}}} y_k(\hat{\boldsymbol x})^{\top}\overline {\boldsymbol d}_i \sum_{i=m_1+1}^{m} \left ( \alpha_i \overline \omega_i + (1-\alpha_i)  \underline \omega_i - \frac{\overline \omega_i + \underline \omega_i}{2}\right)\nonumber\\
		\Leftrightarrow & {Z^*}(m)-	{Z^*}(m_1) \leq  \max_{k=1}^K  \ \max\limits_{\substack{0 \le \alpha_i \le 1\\  \forall i \in \{m_1+1,\dots,m\}}} y_k(\hat{\boldsymbol x})^{\top}\overline {\boldsymbol d}_i \sum_{i=m_1+1}^{m} \left ( (\alpha_i - \frac{1}{2})\overline \omega_i + (\frac{1}{2}-\alpha_i)  \underline \omega_i \right)\label{sub1}\\
		\Leftrightarrow & {Z^*}(m)-	{Z^*}(m_1) \leq  \max_{k=1}^K  \ \sum_{i=m_1+1}^{m} \max \left\{ y_k(\hat{\boldsymbol x})^{\top}\overline {\boldsymbol d}_i  \left ( - \frac{1}{2}\overline \omega_i + \frac{1}{2} \underline \omega_i \right), y_k(\hat{\boldsymbol x})^{\top}\overline {\boldsymbol d}_i  \left (\frac{1}{2}\overline \omega_i -\frac{1}{2}\underline \omega_i \right)\right\}\label{sub2}\\
		\Leftrightarrow & {Z^*}(m)-	{Z^*}(m_1) \leq  \max_{k=1}^K  \ \sum_{i=m_1+1}^{m}  |y_k(\hat{\boldsymbol x})^{\top}\overline {\boldsymbol d}_i|  \left (\frac{\overline \omega_i- \underline \omega_i}{2}\right).\nonumber
	\end{align}
\end{small}
Note that \eqref{sub2} is equivalent to \eqref{sub1} because the inner maximization problem in \eqref{sub1} is linear. Therefore, its optimal solution is one of the extreme points of $0\leq \alpha_i \leq 1$, i.e., either $\alpha_i=0$ or $\alpha_i=1$.\endproof

The theoretical upper bound developed in Theorem \ref{ROLBTheoretical0} brings two benefits: (i) it provides a rough approximation for the optimal value of Problem \eqref{RO1}, which may not be solved efficiently in practice; and (ii) it determines how many principal components are required to reach a preferred gap, demonstrating a trade-off between computational burden and solution quality.

\section{Probabilistic Guarantees}\label{Section 4}

In this section, we derive probabilistic guarantees for the performance of $\mathcal U_{\textsc{pca}}(\mathcal S, m_1)$ when all the principal components are utilized to construct these uncertainty sets, i.e., $m_1 = m$. To that end, we develop explicit lower bounds on the number of scenario samples required to construct these sets with desired probabilistic performance. To derive the probabilistic guarantees, we consider  no assumptions on the probability distribution of uncertainty data.

\begin{theorem}\label{third-thm-3}
	If $N \geq N^*(m)=\lceil\frac{1}{\epsilon}\frac{e}{e-1}(2m-1+\ln{\frac{1}{\beta}})\rceil$, then we have $1-\beta$ confidence that any realization $\boldsymbol s$ belongs to the uncertainty set $\ \mathcal U_{\textsc{pca}}(\mathcal S, m)$ with the probability of at least $1-\epsilon$, i.e.,
	\begin{equation}\label{ineq}
		\Prob_{\hat{\boldsymbol{s}}}\{\Prob_{\boldsymbol s}\{\boldsymbol s \in \mathcal U_{\textsc{pca}}(\mathcal S, m)\}\geq 1-\epsilon\}\geq 1-\beta,\nonumber
	\end{equation}
	where $0< \epsilon <1$, $0< \beta <1$, and $\hat{\boldsymbol{s}}=\{\boldsymbol s_1,\cdots,\boldsymbol s_N\}$. 	
\end{theorem}
\proof{Proof.}
The result is deduced from \cite{margellos2014road} and thus we omit the proof. \endproof
From Theorem~\ref{third-thm-3}, when $m=1$, $N^*(m)$ becomes $N^*(1)=\lceil\frac{1}{\epsilon}\frac{e}{e-1}(1+\ln{\frac{1}{\beta}})\rceil$.

\begin{theorem}\label{third-thm-4}
When $m=1$. If $N\geq \frac{\ln{{\beta}}-\ln{{(1-\epsilon+N^*(1)\epsilon)}}}{\ln{(1-\epsilon)}}+1$, then we have $1-\beta$ confidence that any realization $\boldsymbol s$ belongs to the uncertainty set $\ \mathcal U_{\textsc{pca}}(\mathcal S,m)$ with the probability of at least $1-\epsilon$, i.e.,
	\begin{equation}\label{ROineq}
		\Prob_{\hat{\boldsymbol{s}}}\{\Prob_{\boldsymbol s}\{\boldsymbol s \in \mathcal U_{\textsc{pca}}(\mathcal S, m)\}\geq 1-\epsilon\}\geq 1-\beta,
	\end{equation}
	where $0< \epsilon <1$, $0< \beta <1$, and $\hat{\boldsymbol{s}}=\{\boldsymbol s_1,\cdots,\boldsymbol s_N\}$. 	
\end{theorem}
\proof{Proof.}

With $m=1$, we have 
\begin{align}
	\Prob_{\boldsymbol s}\{\boldsymbol s \in \mathcal U_{\textsc{pca}}(\mathcal S, 1)\} &= \Prob_{\boldsymbol s}\{\boldsymbol{\bar s}_j\leq \boldsymbol{\bar s}\leq \boldsymbol{\bar s}_j,\quad \forall j \in [N]\} \nonumber\\
	& = \Prob_{\boldsymbol s}\{\min_{j=1}^{N}\boldsymbol{\bar s}_j \leq \boldsymbol{\bar s}\leq \max_{j=1}^{N}\boldsymbol{\bar s}_j\}\nonumber
\end{align}
Thus we have a sufficient and necessary condition for $\Prob_{\boldsymbol s}\{\boldsymbol s \in \mathcal U_{\textsc{pca}}(\mathcal S, 1)\}\geq 1-\epsilon$: 
$$\Prob_{\boldsymbol s}\{\min_{j=1}^{N}\boldsymbol{\bar s}_j \leq \boldsymbol{\bar s}\leq \max_{j=1}^{N}\boldsymbol{\bar s}_j\} \geq 1 - \epsilon.$$
Similarly, we have
\begin{scriptsize}
	\begin{align}
		\Prob_{\hat{\boldsymbol{s}}}\left\{\Prob_{\boldsymbol s}\{\min_{j=1}^{N}\boldsymbol{\bar s}_j \leq \boldsymbol{\bar s}\leq \max_{j=1}^{N}\boldsymbol{\bar s}_j\} \geq 1 - \epsilon\right\} &= \Prob_{\hat{\boldsymbol{s}}}\left\{F_1( \max_{j=1}^{N} \boldsymbol{\bar s}_j)-F_1( \min_{j=1}^{N} \boldsymbol{\bar s}_j)\geq 1 - \epsilon\right\} \nonumber\\
		&= \Prob_{\hat{\boldsymbol{s}}}\left\{\max_{j=1}^{N} F_1( \boldsymbol{\bar s}_j)-\min_{j=1}^{N} F_1(\boldsymbol{\bar s}_j) \geq 1 - \epsilon,\right\} \nonumber	
	\end{align}
\end{scriptsize}
where the second equality is due to the fact that $F_1$ is non-decreasing. 
Thus, to make inequality~\eqref{ROineq} hold, it is equivalent to have 
$$\Prob_{\hat{\boldsymbol{s}}}\left\{\max_{j=1}^{N}F_1(\boldsymbol{\bar s}_j)-\min_{j=1}^{N}F_1( \boldsymbol{\bar s}_j)\geq 1 - \epsilon\right\} \geq 1-\beta.$$
As $F_1(\boldsymbol{\bar s}_j)$ is a random variable with the probability distribution of $\text{unif}(0,1)$, by defining $\xi_1=\max_{j=1}^{N}F_1(\boldsymbol{\bar s}_j)$ and $\xi_2=\min_{j=1}^{N}F_1( \boldsymbol{\bar s}_j)$, the joint probability density function of order statistics $\xi_1$ and $\xi_2$ is $N(N-1)(\xi_1-\xi_2)^{N-2}$ when $\xi_1\geq \xi_2$ and zero otherwise \citep{casella2021statistical}. Consequently, we have
\begin{align}
	\Prob_{\hat{\boldsymbol{s}}}\left\{\max_{j=1}^{N}F_1(\boldsymbol{\bar s}_j)-\min_{j=1}^{N}F_1( \boldsymbol{\bar s}_j)\geq 1 - \epsilon\right\} &=\int_{0}^{\epsilon}\int_{1-\epsilon+\xi_2}^{1} N(N-1)(\xi_1-\xi_2)^{N-2}d\xi_1d\xi_2\nonumber\\
	& =1-\left(1-\epsilon\right)^{N}-N\left(1-\epsilon\right)^{N-1}\epsilon.\nonumber
\end{align}
Moreover, $1-\left(1-\epsilon\right)^{N}-N\left(1-\epsilon\right)^{N-1}\epsilon\nonumber  \geq 1-\beta$ is equivalent to 
\begin{equation}\label{third-paper-ineq}
	\left(1-\epsilon\right)^{N-1}(1-\epsilon+N\epsilon)  \leq \beta.
\end{equation}
Thus, to complete the proof, it is sufficient to show $N =\frac{\ln{{\beta}}-\ln{{(1-\epsilon+N^*(1)\epsilon)}}}{\ln{(1-\epsilon)}}+1$ satisfies the inequality \eqref{third-paper-ineq}. Based on the results of Theorem~\ref{third-thm-3}, we know that $N^*(1)$ satisfies the inequality \eqref{third-paper-ineq}. That means there exists an $\bar N \leq N^*(1)$, such that the inequality \eqref{third-paper-ineq} holds with $N=\bar N$, i.e.,
\begin{equation}\label{third-paper-ineq-1}
	\left(1-\epsilon\right)^{\bar N-1}(1-\epsilon+\bar N\epsilon)  \leq \beta,
\end{equation}
The above inequality \eqref{third-paper-ineq-1} can be implied by the condition  $\left(1-\epsilon\right)^{\bar N-1}(1-\epsilon+N^*(1)\epsilon)  \leq \beta$, which is equivalent to $\bar N\geq \frac{\ln{{\beta}}-\ln{{(1-\epsilon+N^*(1)\epsilon)}}}{\ln{(1-\epsilon)}}+1$. Therefore the proof is complete.

\endproof

It is worth noting that the derived lower bound $N_1^{*}=\lceil\frac{\ln{{\beta}}-\ln{{(1-\epsilon+N^*(1)\epsilon)}}}{\ln{(1-\epsilon)}}+1\rceil$ for $N$ in Theorem~\ref{third-thm-4} is always smaller than the existing lower bound $N^{*}(1)=\lceil \frac{1}{\epsilon}\frac{e}{e-1}(1+\ln{\frac{1}{\beta}}) \rceil$ in Theorem \ref{third-thm-3} when $m=1$. For instance, when $\alpha=0.1$, $\beta=0.1$ and $m=1$, $N_1^{*}=41$ is while $N_1^{*}(1)$ is $53$. Therefore, the developed lower bound improves the existing work. We next extend the results of Theorem~\ref{third-thm-4}  to the general case of $m$. 

\begin{corollary}\label{RO-Corrolary1}
	Let $N^*= \lceil \frac{m}{\epsilon}\frac{e}{e-1}(1+\ln{\frac{m}{\beta}})\rceil$. If $N\geq \frac{\ln{{\beta}}-\ln{{(m-\epsilon+N^*\epsilon)}}}{\ln{(1-\frac{\epsilon}{m})}}+1$, then we have $1-\beta$ confidence that any realization $\boldsymbol s$ belongs to the uncertainty set $\ \mathcal U_{\textsc{pca}}(\mathcal S, m)$ with the probability of at least $1-\epsilon$, i.e.,
	\begin{equation}\label{ROineq1}
		\Prob_{\hat{\boldsymbol{s}}}\{\Prob_{\boldsymbol s}\{\boldsymbol s \in \mathcal U_{\textsc{pca}}(\mathcal S, m)\}\geq 1-\epsilon\}\geq 1-\beta,
	\end{equation}
	where $0< \epsilon <1$, $0< \beta <1$, and $\hat{\boldsymbol{s}}=\{\boldsymbol s_1,\cdots,\boldsymbol s_N\}$. 	
\end{corollary}
\proof{Proof.}

We have 
\begin{align}
	\Prob_{\boldsymbol s}\{\boldsymbol s \in \mathcal U_{\textsc{pca}}(\mathcal S, m)\} &= \Prob_{\boldsymbol s}\{\boldsymbol{\bar s}_j^i \leq \boldsymbol{\bar s}^i\leq \boldsymbol{\bar s}_j^i,\quad \forall i \in [m], \ \forall j \in [N]\} \nonumber\\
	& = \Prob_{\boldsymbol s}\{\min_{j=1}^{N}\boldsymbol{\bar s}_j^i \leq \boldsymbol{\bar s}^i\leq \max_{j=1}^{N}\boldsymbol{\bar s}_j^i,\quad \forall i \in [m]\}\nonumber\\
	& \geq \sum_{i=1}^m \Prob_{\boldsymbol s}\{\min_{j=1}^{N}\boldsymbol{\bar s}_j^i \leq \boldsymbol{\bar s}^i\leq \max_{j=1}^{N}\boldsymbol{\bar s}_j^i\}-m+1\nonumber
\end{align}
where the last inequality is due to the Bonferroni Inequalities. Thus we have a sufficient condition for $\Prob_{\boldsymbol s}\{\boldsymbol s \in \mathcal U_{\textsc{pca}}(\mathcal S, m)\}\geq 1-\epsilon$: 
$$\sum_{i=1}^m \Prob_{\boldsymbol s}\{\min_{j=1}^{N}\boldsymbol{\bar s}_j^i \leq \boldsymbol{\bar s}^i\leq \max_{j=1}^{N}\boldsymbol{\bar s}_j^i\}-m+1  \geq 1 - \epsilon,$$
which can be further implied by: 
$$\Prob_{\boldsymbol s}\{\min_{j=1}^{N}\boldsymbol{\bar s}_j^i \leq \boldsymbol{\bar s}^i\leq \max_{j=1}^{N}\boldsymbol{\bar s}_j^i\} \geq 1 - \frac{\epsilon}{m},\quad \forall i \in [m]. \nonumber$$
Similarly, we have
\begin{scriptsize}
	\begin{align}
		\Prob_{\hat{\boldsymbol{s}}}\left\{\Prob_{\boldsymbol s}\{\min_{j=1}^{N}\boldsymbol{\bar s}_j^i \leq \boldsymbol{\bar s}^i\leq \max_{j=1}^{N}\boldsymbol{\bar s}_j^i\} \geq 1 - \frac{\epsilon}{m}, \ \forall i \in [m]\right\} &= \Prob_{\hat{\boldsymbol{s}}}\left\{F_i( \max_{j=1}^{N} \boldsymbol{\bar s}_j^i)-F_i( \min_{j=1}^{N} \boldsymbol{\bar s}_j^i)\geq 1 - \frac{\epsilon}{m}, \ \forall i \in [m]\right\} \nonumber\\
		& \geq \sum_{i=1}^m \Prob_{\hat{\boldsymbol{s}}}\left\{\max_{j=1}^{N} F_i( \boldsymbol{\bar s}_j^i)-\min_{j=1}^{N} F_i(\boldsymbol{\bar s}_j^i) \geq 1 - \frac{\epsilon}{m}\right\} - m + 1\nonumber.			
	\end{align}
\end{scriptsize}
Thus, to make inequality~\eqref{ROineq1} hold, it is sufficient to have 
$$\Prob_{\hat{\boldsymbol{s}}}\left\{\max_{j=1}^{N}F_i(\boldsymbol{\bar s}_j^i)-\min_{j=1}^{N}F_i( \boldsymbol{\bar s}_j^i)\geq 1 - \frac{\epsilon}{m}\right\} \geq 1-\frac{\beta}{m}.$$
Then by the results of Theorem~\ref{third-thm-4}, the conclusion follows. 
\endproof


It should be noted that, in contrast to the case $m=1$,  the derived lower bound $N_1^{**}=\lceil\frac{\ln{{\beta}}-\ln{{(m-\epsilon+N^*\epsilon)}}}{\ln{(1-\frac{\epsilon}{m})}}+1\rceil$ for $N$ in Corollary~\ref{RO-Corrolary1} is not always smaller than the existing lower bound $N_2^{**}=\lceil \frac{1}{\epsilon}\frac{e}{e-1}(2m-1+\ln{\frac{1}{\beta}}) \rceil$ in Theorem \ref{third-thm-3} when $m \geq 2$. Therefore, the developed lower bound complements the existing work.

\section{Computational Experiments}\label{ROExperiments}
We conduct comprehensive computational experiments to show the effectiveness of the proposed scenario-induced uncertainty sets using two applications: Knapsack and power grid problems. The mathematical models are implemented by MATLAB R2021a (ver. 9.10) API of 
Gurobi  (ver. 9.1) on a PC with a 64-bit Windows Operating System, an Intel(R) Core(TM) i7-7700 CPU @
3.60 GHz processor, and 16 GB RAM. In Section
\ref{Computational Setup}, we specify the proposed uncertainty sets in the context of the Knapsack and power grid problems. In Section \ref{Computational Results}, we address how to randomly generate test instances of these applications and  report the numerical results along with their analyses.

\subsection{Computational Setup}\label{Computational Setup}

In this section, we specify the proposed $\mathcal U_{\textsc{pca}}(\mathcal S, m_1)$, and  $\mathcal U_\cap (\mathcal S, m_1)$ uncertainty sets in the context of the knapsack and power grid problems. 
\subsubsection{Knapsack Problem}\label{knapsack problem}
We are given a set of items, each with a given value and uncertain weight, that we wish to pack into a container with a maximum capacity limit. The goal is to maximize the total value of the packed items by choosing a subset of the items that fit into the container.
\begin{align}\label{Knapsack 1}
	\max & \quad \boldsymbol v^\top \boldsymbol x \\	
	\text{s.t.} &\quad \boldsymbol w^\top  \boldsymbol x \leq W, \quad \forall \boldsymbol w \in \mathcal U, \nonumber\\	
	& \quad x_z\in \{0,1\} , \quad \forall z \in[n].\nonumber
\end{align} 
This problem is a static RO problem. In Problem \eqref{Knapsack 1}, parameter $n$ represents the number of items and $\boldsymbol v \in \R^n$ and $\boldsymbol w \in \R^n$ denote the values and weights of the items, respectively. Parameter $W$ indicates the maximum capacity of the container and $\mathcal U$ represents the uncertainty set of the uncertain weights. Decision variable $x_z, \ \forall z \in [n]$, indicates if item $z$ is packed into the container (i.e., $x_z=1$) or not (i.e., $x_z=0$). The objective is to maximize the total value of the packed items subject to the constraint that the total weight of the packed items does not exceed the maximum capacity of the container. Problem \eqref{Knapsack 1} can be reformulated as the following bi-level problem:   
\begin{subequations}\label{Knapsack 2}
	\begin{align}
		\max & \quad \boldsymbol v^\top \boldsymbol x \\	
		\text{s.t.} &\quad \max_{\boldsymbol w\in \mathcal U} \quad \boldsymbol w^\top  \boldsymbol x \leq W,\label{Knapsack 2:Const 1}\\	
		& \quad x_z\in \{0,1\} , \quad \forall z\in[n].\nonumber
	\end{align}
\end{subequations}
After applying $\mathcal U_{\textsc{pca}}(\mathcal S, m_1)$ set to Problem \eqref{Knapsack 2} and replacing the inner optimization problem \eqref{Knapsack 2:Const 1} with its dual formulation, Problem \eqref{Knapsack 2} is equivalent to the following problem:
\begin{align}\label{RKP-ASIU}
	\max  \quad & \boldsymbol v^\top \boldsymbol x \nonumber\\
	\text{s.t.} \quad & \sum_{i=1}^{m_1} \beta_i + \sum_{i=1}^{m_1} \frac{\underline \omega_i}{||\boldsymbol d_i||}\boldsymbol d_i^\top\x +\sum_{i=m_1+1}^{m} \frac{\overline \omega_i + \underline \omega_i}{2||\boldsymbol d_i||}\boldsymbol d_i^\top\x\leq W \nonumber\\
	& \beta_i\geq \frac{\overline \omega_i-\underline \omega_i}{||\boldsymbol d_i||}\boldsymbol d_i^\top\x, \quad \forall i\in[m_1],  \nonumber\\
	&\beta_i \in \R_+, \quad \forall i\in[m_1],  \nonumber\\
	& x_z\in \{0,1\} , \quad \forall z\in[n], \nonumber
\end{align}
where $\beta_i, \ \forall i \in [m_1]$ are dual decision variables.  Moreover, after applying uncertainty set $\mathcal U_\cap(\mathcal S, m_1)$ to Problem \eqref{Knapsack 2} and replacing the inner optimization problem \eqref{Knapsack 2:Const 1} with its dual formulation, Problem \eqref{Knapsack 2} is reformulated as the following problem:
\begin{align}
	\max  \quad & \boldsymbol v^\top \boldsymbol x \nonumber\\
	\text{s.t.} \quad & \sum_{i=1}^{m_1} \beta_i + (\boldsymbol U-\boldsymbol c)^\top \boldsymbol \gamma+(\boldsymbol c-\boldsymbol L)^\top \boldsymbol \zeta + \boldsymbol c_0^\top\x \leq W, \nonumber\\
	& \beta_i + \frac{\overline \omega_i-\underline \omega_i}{||\boldsymbol d_i||}\boldsymbol d_i^\top(\boldsymbol \gamma-\boldsymbol \zeta) \geq \frac{\overline \omega_i-\underline \omega_i}{||\boldsymbol d_i||}\boldsymbol d_i^\top\x, \quad  \forall i\in[m_1], \nonumber\\
	& \boldsymbol \gamma,\boldsymbol \zeta \in \R^{n}_+,\nonumber\\
	& \beta_i \in \R_+, \ \forall i\in[m_1], \nonumber\\
	& x_z\in \{0,1\}, \  \forall z\in[n], \nonumber
\end{align}
where $\boldsymbol U$ and $\boldsymbol L$ are respectively the upper bound and lower bound of the box uncertainty set, 
\[\boldsymbol c=\sum_{i=1}^{m_1} \frac{\underline \omega_i }{||\boldsymbol d_i||}\boldsymbol d_i+\sum_{i=m_1+1}^{m} \frac{\overline \omega_i + \underline \omega_i}{2||\boldsymbol d_i||}\boldsymbol d_i, \]
$\beta_i, \ \forall i \in [m_1]$, $\boldsymbol \gamma$, and $\boldsymbol \zeta$  are dual decision variables.

\subsubsection{Power Grid Problem}\label{Power problem}
In this problem, we consider a dispatchable power grid. This power grid is a network of generator stations, transmission systems, and consumers that delivers power from generators to consumers. The power generated by generators is referred to as output while consumers' power demand is referred to as load. Load is also considered as any component of the power grid that consumes power. A bus is defined as a vertical line at which several components of a power grid such as loads or generators are connected. This power grid includes load buses and generator buses. The generators of this power grid are dispatchable, i.e., they can be dispatched on demand by adjusting their output according to power orders. Moreover, the dispatchable generators are subject to ramping constraints.  A ramp event is defined as a power increase or decrease event that happens in a time unit. More specifically, a ramp-up event occurs when power increases while a ramp-down event occurs when the power decreases.

To balance output and load, the load shedding and output curtailment procedures are performed on this power grid. Load shedding is the act of switching off output to some consumers when the load is more than output to prevent the power grid from collapsing. Output curtailment is the act of deliberately reducing output to below what could have been generated due to the maintenance of the transmission system or the overloaded transmission system when output is more than load. We assume there is no power flow limitation for this power grid. With this background, the power grid problem with uncertain load (demand) is defined as follows:
\begin{subequations}\label{PGA-ori}
	\begin{align}
		\max_{\hat d \in \mathcal U} \ \min_{\hat p, \hat {\overline q}, \hat {\underline q}}  \quad &  \sum_{t=1}^T\sum_{g=1}^{G} c_g p_g^t+ \overline M\sum_{t=1}^T\overline q^t+  \underline M \sum_{t=1}^T\underline q^t\label{PGA-ori:obj}\\
		\text{s.t.} \quad & \sum_{g=1}^{G} p_g^t + \overline q^t - \underline q^t = \sum_{l=1}^{L} d_l^t, \quad  \forall t\in[T], \label{PGA-ori:const1} \\
		&- \underline R_g \le p_g^t - p_g^{t-1} \le \overline R_g, \quad \forall g \in [G], \  \forall t\in [T],\label{PGA-ori:const2}\\
		& p_g^t \le \overline P_g, \quad \forall g  \in [G], \  \forall t\in[T],\label{PGA-ori:const3}\\
		& p_g^t, \overline q^t, \underline q^t \in \R_+, \quad \forall g \in [G], \ \forall t \in [T].\nonumber 
	\end{align}
\end{subequations}
where $\hat d = \{d_l^t, \ \forall l \in [L], \forall t \in [T]\}$ $\hat p = \{p_g^t, \ \forall g \in [G], \forall t \in [T]\}$, $\hat {\overline q} = \{\overline q^t, \ \forall t \in [T]\}$, and $\hat {\underline q} = \{\underline q^t, \ \forall t \in [T]\}$ . This problem is a special case of an ARO problem. In Problem \eqref{PGA-ori}, parameters $T$, $G$, $L$ indicate the total number of time units, generator buses, and load buses while each time unit, generator bus, and load bus is identified by indices $t$, $g$, and $l$, respectively. The cost of generating one megawatt of output by generator bus $g$ in one time unit is denoted by $c_g$ and parameters $\overline M$ and $\underline M$ indicate the penalty costs for one megawatt of load shedding and output curtailment performed in one time unit, respectively. The uncertain load (demand) of load bus $l$ at time $t$ is represented by $d_l^t$.  Parameters $\underline  R_g$ and $\overline R_g$ denote the maximum allowed ramp-down and ramp-up in megawatt between two consecutive time units for generator bus $g$ and parameter $\overline P_g$ represents the maximum output capacity of generator bus $g$ in each time unit.  Decision variable $p_g^t$ denotes the generated output by generator bus $g$ at time $t$ in megawatt. Decision variables $\overline q^t$ and $\underline q^t$ represent the performed load shedding and output curtailment on the power grid at time $t$ in megawatt. 

Problem \eqref{PGA-ori} minimizes the worst-case total economic dispatch cost, including output generation cost,  load shedding penalty cost, and output curtailment penalty cost, by determining the optimal generated output in megawatt by each generator bus in each time unit and optimal performed output curtailment and load shedding in megawatt in each time unit. Constraint \eqref{PGA-ori:const1} balances the total output of generator buses and the total load of load buses considering performed output curtailment and load shedding in each time unit. Constraint \eqref{PGA-ori:const2} limits the ramp of each generator bus in each two consecutive time units to its maximum allowed ramp-up and ramp-down. Constraint \eqref{PGA-ori:const3} guarantees that the generated output by each generator bus in each time unit does not exceed its maximum output generation capacity.

After replacing the inner optimization problem with its dual formulation, Problem \eqref{PGA-ori} is equivalent to the following problem:
\begin{align}\label{power grid 2}
	\max_{\hat d \in \mathcal U} \ \max_{\hat{\overline{y}},\hat{\underline{y}}, \hat z, \hat x}  \quad &   \sum_{t=1}^T\left(\sum_{l=1}^{L} d_l^t\right)x^t - \sum_{t=1}^T\sum_{g=1}^{G}(\overline R_g \overline{y}_g^t+\underline R_g \underline{y}_g^t+\overline P_g z_g^t)   \\
	\text{s.t.} \quad & x^t- \overline{y}_g^t+\overline{y}_g^{t+1}+\underline{y}_g^t-\underline{y}_g^{t+1}-z_g^t \leq c_g, \quad \forall g \in [G], \  \forall t\in[T-1], \nonumber\\
	&  x^T- \overline{y}_g^T+\underline{y}_g^T-z_g^T\leq c_g, \quad \forall g \in [G], \nonumber\\
	& -\underline M \leq x^t \le \overline M, \quad \forall t\in[T],\nonumber\\
	& \overline{y}_g^t, \underline{y}_g^t, z_g^t \in \R_+, \quad \forall g \in [G], \  \forall t \in [T], \nonumber
\end{align}
where $\hat{\overline{y}} = \{\overline{y}_g^t, \ \forall g \in [G], t \in [T]\}$, $\hat{\underline{y}}= \{\underline{y}_g^t, \ \forall g \in [G], t \in [T]\}$, $\hat z= \{z_g^t, \ \forall g \in [G], t \in [T]\}$, and $\hat x= \{x^t, \ \forall t \in [T]\}$ are dual decision variables. Then, applying $\mathcal U_{\textsc{pca}}(\mathcal S, m_1)$ set to Problem \eqref{power grid 2} leads to the followig problem:
\begin{align}\label{PGA-ASIU}
	\max_{\hat{\overline{y}},\hat{\underline{y}}, \hat z, \hat x, \hat{\alpha}}  \quad &  \sum_{t=1}^T\left(\sum_{l=1}^{L}\boldsymbol u _{(t-1)L+l}\right)x^t - \sum_{t=1}^T\sum_{g=1}^{G}(\overline R_g \overline{y}_g^t+\underline R_g \underline{y}_g^t+\overline P_g z_g^t)\\
	\text{s.t.} \quad & x^t- \overline{y}_g^t+\overline{y}_g^{t+1}+\underline{y}_g^t-\underline{y}_g^{t+1}-z_g^t \leq c_g, \quad \forall g \in [G], \  \forall t\in[T-1], \nonumber\\
	& x^T- \overline{y}_g^T+\underline{y}_g^T-z_g^T\leq c_g, \quad \forall g \in [G], \nonumber\\
	& -\underline M \leq x^t \le M, \quad  \forall t\in[T],\nonumber\\
	& \overline{y}_g^t, \underline{y}_g^t, z_g^t \in \R_+, \quad \forall g \in [G], \  \forall t \in [T], \nonumber\\
	&0\leq \alpha_i\leq 1, \quad \forall i\in[m_1], \nonumber
\end{align}
where $\hat \alpha = \{\alpha_i, \ \forall i \in [m_1]\}$ and 
\[ \uu = \sum_{i=1}^{m_1} \ \left(\alpha_i \left(\frac{\overline \omega_i }{||\boldsymbol d_i||}\boldsymbol d_i\right) + (1-\alpha_i)\left( \frac{ \underline \omega_i}{||\boldsymbol d_i||}\boldsymbol d_i\right)\right) + \sum_{i=m_1+1}^{m} \frac{\overline \omega_i + \underline \omega_i}{2||\boldsymbol d_i||}\boldsymbol d_i. \] 
Uncertain parameters  $d_l^t, \ \forall l \in [L], \forall t \in [T]$, can be represented alternatively as a single uncertain parameter (random variable) vector $\boldsymbol {u} \in \R^m$, where $m=TL$. Accordingly, uncertain parameter $d_l^t$ for a specific $l$ and $t$ is located in the $((t-1)L+l)^{th}$ element of vector $\uu \in \R^{TL}$. Therefore, $\boldsymbol u_{(t-1)L+l}$ in the objective function denotes the $((t-1)L+l)^{th}$ element of $\boldsymbol u$, which is $d_l^t$ for a given $t$ and $l$.  

When $\hat{\overline{y}},\hat{\underline{y}}, \hat z$, and $\hat x$ are fixed, Problem \eqref{PGA-ASIU} is equivalent to the following problem:
\begin{align}\label{PGA-ASIU2}
	\max_{\hat \alpha }  \quad &  \sum_{t=1}^T\left(\sum_{l=1}^{L}\boldsymbol u_{(t-1)L+l}\right)x^t - \sum_{t=1}^T\sum_{g=1 }^{G}(\overline R_g \overline{y}_g^t+\underline R_g \underline{y}_g^t+\overline P_g z_g^t)  \nonumber \\
	\text{s.t.} \quad & x^t- \overline{y}_g^t+\overline{y}_g^{t+1}+\underline{y}_g^t-\underline{y}_g^{t+1}-z_g^t \leq c_g, \quad \forall g \in [G], \ \forall t\in[T-1],  \nonumber\\
	& x^T- \overline{y}_g^T+\underline{y}_g^T-z_g^T\leq c_g, \quad \forall g \in [G],  \nonumber\\
	&  -\underline M \leq  x^t \le \overline M, \quad   \forall t\in[T],  \nonumber\\
	& \overline{y}_g^t, \underline{y}_g^t, z_g^t \in \R_+, \quad \forall g \in [G], \  \forall t \in [T], \nonumber\\
	&\alpha_i\in\{0,1\}, \quad \forall i\in[m_1], \nonumber
\end{align}
where $\hat \alpha$ is the only decision variable of this linear problem and its feasible region is a box. Therefore,  Problem \eqref{PGA-ASIU}, which is bi-level, is equivalent to the following single-level mixed integer linear problem:
\begin{align}
	\max_{\hat{\overline{y}},\hat{\underline{y}}, \hat z, \hat v, \hat x, \hat{\alpha}}  \quad &  \sum_{t=1}^T\sum_{l=1}^{L}\boldsymbol q_{(t-1)L+l} - \sum_{t=1}^T\sum_{g=1}^{G}(\overline R_g \overline{y}_g^t+\underline R_g \underline{y}_g^t+\overline P_g z_g^t)  \nonumber \\
	\text{s.t.} \quad & x^t- \overline{y}_g^t+\overline{y}_g^{t+1}+\underline{y}_g^t-\underline{y}_g^{t+1}-z_g^t \leq c_g, \quad \forall g \in [G],\  \forall t\in[T-1], \nonumber\\
	& x^T- \overline{y}_g^T+\underline{y}_g^T-z_g^T\leq c_g, \quad \forall g \in [G], \nonumber\\
	&  -\underline M \leq x^t \le \overline M,\quad   \forall t\in[T], \nonumber\\
	& -\underline{M}\alpha_i\leq  v_i^t \leq \overline{ M}\alpha_i,\quad \forall i\in[m_1], \  \forall t\in[T], \nonumber \\
	& x^t-(1-\alpha_i)\overline{M}\leq v_i^t\leq x^t+(1-\alpha_i)\underline{M},\quad \forall i\in[m_1], \  \forall t\in[T], \nonumber\\
	& \overline{y}_g^t, \underline{y}_g^t, z_g^t, v_i^t \in \R_+, \quad \forall g \in [G], \  \forall t \in [T], \ \forall i \in [m_1], \nonumber\\
	&\alpha_i\in\{0,1\},\quad \forall i\in[m_1], \nonumber
\end{align}
where $\hat v = \{v^t_i, \ \forall i \in [m_1], \forall t \in [T]\}$ and
\[ \boldsymbol q = \sum_{i=1}^{m_1} \left ( v_i^t\left(\frac{\overline \omega_i}{||\boldsymbol d_i||}\boldsymbol d_i\right) + (x^t-v_i^t)\left( \frac{\underline \omega_i}{||\boldsymbol d_i||}\boldsymbol d_i\right)\right) + x^t \left(\sum_{i=m_1+1}^{m} \frac{\overline \omega_i + \underline \omega_i}{2||\boldsymbol d_i||}\boldsymbol d_i\right).\]

\subsection{Computational Results}\label{Computational Results}

We first explain how we generated random test instances for the Knapsack and power grid problems used to evaluate the performance of the proposed SIU-based RO models. Then, we compare the robust counterparts using these uncertainty sets in terms of the conservativeness of their solutions and the computational time needed to solve them to optimality. Finally, we further investigate the performance of the developed uncertainty sets by performing sensitivity analysis concerning the parameters of the uncertainty sets and the parameters of the Knapsack and power grid problems.

\subsubsection{Instance Generation and Table Header Description} We conduct our computational experiments
to solve various instances of the robust knapsack and power grid problems with the proposed uncertainty sets.
To generate test instances of the knapsack problem, we follow the same experimental setup proposed by \cite{bertsimas2004price}. As this problem is NP-hard, we generate random Knapsack problems of size $n = 200$, which can be solved to optimality by
off-the-shelf optimization solvers. The value of each item, i.e., $v_z, \ \forall z \in [200]$, is 
randomly selected from the set $\{16,17,\dots,77\}$. The weight of each item, i.e., $w_z, \ \forall z \in [200]$, is assumed to be uncertain, dependent on the weights of other items, and follows a Normal distribution, i.e., $w_z \sim N(\mu_z, \sigma^2_z), \ \forall z \in [200]$, where $\mu_z$ and $\sigma^2_z$ denote its mean and variance, respectively. Parameter $\mu_z$ is randomly chosen from the set $\{20, 21, \dots, 29\}$ and  $\sigma^2_z$ is assumed to be  $\sigma^2_z=\frac{\mu_z^2}{300}$. There exists the same correlation between any two dependent weights, which is denoted by $\rho$ and implies their dependency. More specifically, the weight of item $z$ is assumed to be correlated with the weight of item $z+1$ for any odd values of $z$. With this assumption and the definition of correlation, i.e., $\rho = \frac{\text{Cov}(w_z,w_{z+1})}{\sigma_z\sigma_{z+1}}$, the covariance matrix of the weights is as follows: 
\begin{align*}
	\text{Cov}(\boldsymbol w)= \left[\begin{smallmatrix} \frac{\mu_1^2}{300}&\frac{\rho\mu_1\mu_2}{300}&0&\dots& 0& 0\\
		\frac{\rho\mu_1\mu_2}{300}& \frac{\mu_2^2}{300} &0&\dots & 0& 0\\
		0& 0& & & 0& 0\\
		\vdots& \vdots &  & \ddots\ddots & \vdots & \vdots\\
		0& 0& & & 0& 0\\
		0 & 0  & \dots &0&\frac{\mu_{199}^2}{300} &\frac{\rho\mu_{199}\mu_{200}}{300}\\
		0 & 0  & \dots &0 &\frac{\rho\mu_{199}\mu_{200}}{300}&\frac{\mu_{200}^2}{300} \\
	\end{smallmatrix}\right].
\end{align*}
To randomly generate correlated normally distributed scenarios, we first generate uncorrelated scenarios, denoted by $\boldsymbol s_j^{''} \in \R^{200}, \forall j \in [N]$, using $w_z \sim N(\mu_z, \sigma^2_z), \ \forall z \in [200]$. Then, we obtain matrix $\boldsymbol M$ by the Cholesky decomposition of $\text{Cov}(\boldsymbol w)$ so that $\boldsymbol M \boldsymbol M^{\top}=\text{Cov}(\boldsymbol w)$. Finally, we generate the correlated scenarios by $\boldsymbol s_j = \boldsymbol \mu + \boldsymbol M \boldsymbol s_j^{''}, \forall j \in [N]$. In the next step, we perform EVD on the sample covariance matrix $\boldsymbol C = \frac{1}{N-1} \boldsymbol X_0^{\top} \boldsymbol X_0$, which is an approximation of $\text{Cov}(\boldsymbol w)$, to develop the proposed uncertainty sets using the PCA technique.

The solutions of the RO Knapsack problems with each proposed uncertainty set are functions of input scenarios. Accordingly, conducting computational experiments based on only one set of scenarios might lead to a biased analysis of the uncertainty set performance. Therefore, we create 10 sets of scenarios, i.e., $\mathcal S_k, \ \forall k \in [10]$, so that each $\mathcal S_k$ includes  $N=10,000$ scenarios of $w_z, \ \forall z \in [200]$, randomly generated  by $N(\mu_z, \sigma^2_z)$. We construct each proposed uncertainty set based on all $\mathcal S_k$ sets, which results in 10 different uncertainty sets of the same type. Then, these uncertainty sets are applied to each Knapsack test problem, leading to 10 problems with the same type of uncertainty set. Finally, the 10 robust knapsack problems are solved to optimality and the average of their optimal objective values and computational times is reported as the performance of the proposed uncertainty set in the context of the given Knapsack test problem.      

To study how the maximum  capacity of the container and the correlation between the weights of items affect the performance of the proposed uncertainty sets, we perform sensitivity analysis with respect to parameters $W$ and $\rho$. To that end,  we conduct our experiments based on three values of the maximum capacity and seven values of the correlation, i.e.,  $W \in \{3000, 4000, 5000\}$ and $\rho \in \{-0.8, -0.5, -0.2, 0, 0.2, 0.5, 0.8\}$.

To generate the test instances of the power grid problem, we consider an IEEE 24-bus system [TK - citation] that consists of 32 generator buses and 17 load buses, planning for a 24-hour horizon, i.e., $T=24$. For simplicity, the 17 load buses are assumed to be grouped into two load buses. Accordingly, the number of generator buses and load buses are set to be 32 and 2, i.e., $G=32$ and $L=2$. For each generator bus, the output generation cost $c_g$ is randomly generated by discrete uniform distribution $\text{unif}(10,150)$ and the maximum output capacity $\overline P_g$ is randomly generated by continuous uniform distribution $\text{unif}(5,245)$. Similarly, the maximum allowed ramp-down and ramp-up for each generator bus, i.e., $\underline  R_g$ and  $\overline R_g$, are randomly generated by continuous uniform distribution $\text{unif}(5,105)$. The penalty costs of the  load shedding and output curtailment are considered $\$500$ and $\$50$ per megawatt, i.e., $\overline M=500$ and $\underline M=50$. 

The two grouped loads in each time, denoted by $d_1^t$ and $d_2^t$, $\forall t \in [24]$, are assumed to be uncertain, dependent on each other, and follow a Normal distribution, i.e., $d_l^t \sim N(\mu_{l}^t, {\sigma_{l}^t}^2), \ \forall l \in [2], \forall t \in [24]$. We assume $\mu_1^t=\frac{2}{5}Q^t$ and $\mu_2^t=\frac{3}{5}Q^t$, where $Q^t$ represents the total loads of 17 load buses in hour $t$ and is randomly generated by continuous uniform distribution $\text{unif}(2100,2900)$. Moreover, the variances of the two grouped loads in each hour are set to be  ${\sigma_{1}^t}^2=\frac{{\mu_1^t}^2}{100}$ and ${\sigma_{2}^t}^2=\frac{{\mu_2^t}^2}{100}$. 

As $d_1^t$ and $d_2^t$, $\forall t \in [24]$, form a collection of random variables or events indexed by different instants of time, we consider them as a stochastic process. This stochastic process is assumed to be a Markov chain, where  each event depends only on the state attained in the previous event. In other words, random variables $(d_1^{t+1}$,$d_2^{t+1})$ only depend on  $(d_1^{t}$,$d_2^{t})$ for all $t \in [23]$. In a Markov chain, there are two types of relationships between random variables, referred to as temporal relationship and spatial relationship. We define the temporal relationship as the correlation between $d_l^{t}$ and $d_l^{t+1}$, which is quantified by the temporal correlation coefficient $\rho_1$ so that $\rho_1 = \frac{\text{Cov}(d_l^{t},d_l^{t+1})}{\sigma_{l}^t\sigma_{l}^{t+1}}, \ \forall l \in [2], \forall t \in [23]$. On the other hand, the spatial relationship is defined as the correlation between $d_1^{t}$ and $d_2^{t}$, which is quantified by the spatial correlation coefficient $\rho_2$ so that $\rho_2 = \frac{\text{Cov}(d_1^{t},d_2^{t})}{\sigma_{1}^t\sigma_{2}^{t}}, \ \forall t \in [24]$.  Uncertain parameters  $d_1^t$ and $d_2^t$, $\forall t \in [24]$, form vector $\boldsymbol {u} \in \R^{48}$ whose $(2(t-1)+l)^{th}$ element is $d_l^t$. Similar to the Knapsack problem, we randomly generate correlated normally distributed scenarios by the randomly generated mean vector, randomly generated uncorrelated scenarios,  and the Cholesky decomposition of $\text{Cov}(\boldsymbol u)$. Then, we perform EVD on the sample covariance matrix $\boldsymbol C$.

To study the effect of temporal and spatial correlations on the performance of the proposed uncertainty sets, we perform sensitivity analysis with respect to parameters $\rho_1$ and $\rho_2$. For this purpose, we consider two and seven settings for the temporal and spatial correlations, respectively , i.e., $\rho_1 \in \{0.5, 0.9\}$ and $\rho_2 \in \{-0.8, -0.5, -0.2, 0, 0.2, 0.5, 0.8\}$. We further investigate the impact of the utilized principal components on the conservativeness of solutions and computational time by setting $m_1$ to $42$ and $36$, i.e., $87.5\%$ and $75\%$ of the $m=48$ principal components.      Like the Knapsack problem, we create 10 sets of $10,000$ scenarios of $d_l^t, \ \forall l \in [2], t \in [24]$, randomly generated by $d_l^t \sim N(\mu_{l}^t, {\sigma_{l}^t}^2)$.

In Section \ref{Table Results}, we will summarize the results of our computational experiments in Tables \ref{table-2} - \ref{table-8}. Columns $``\rho"$, $``\rho_1"$, $``\rho_2"$  report the values of conventional, temporal, and spatial correlations. Column ``Value" represents the optimal objective value of the corresponding RO problem. Column ``Time" reports the computational time of solving the corresponding RO problem to optimality in seconds.  Symbol ``Gap" represents the relative gap in percentage between ``Value'' of the first column and ``Value'' of the corresponding uncertainty set. We define the relative gap between two values as their difference divided by the maximum absolute value. 

\subsubsection{Uncertainty Set Performance}\label{Table Results}

We summarize the results for robust problems with $\mathcal U_{\text{conv}}(\mathcal S)$, $\mathcal U_\text{box}$, $\mathcal U_{\textsc{pca}}(\mathcal S, m)$, and $\mathcal U_\cap (\mathcal S, m)$ uncertainty sets in the context of the Knapsack problem in Tables \ref{table-2} - \ref{table-3}, while Tables \ref{budget-conv-pca-power} - \ref{table-8} report the results for RO problems with $\mathcal U_{\text{conv}}(\mathcal S)$, $\mathcal U_\text{box}$, $\mathcal U_{\textsc{pca}}(\mathcal S, m)$, and $\mathcal U_{\textsc{pca}}(\mathcal S, m_1)$ uncertainty sets in the context of the power grid problem. Sensitivity analysis with respect to the maximum container capacity $W$ and varying values of correlation $\rho$ are reported in Tables \ref{table-2} - \ref{table-3}. Tables \ref{table-4} and \ref{table-8} show how the number of utilized principal components  affects the performance of $\mathcal U_{\textsc{pca}}(\mathcal S, m_1)$, reporting the sensitivity analysis results with respect to $m_1$. 

Sensitivity analyses with respect to the three values of temporal correlation $\rho_1$ are reported through Tables \ref{table-4} - \ref{table-8}. Within each of these tables, we hold $\rho_1$ constant and present sensitivity analysis results for varying values of spatial correlation $\rho_2$. In general, shorter computational times imply greater tractability, and better objective values imply less conservativeness (i.e., a larger value for the Knapsack problem (maximization) and a smaller value for the power grid problem (minimization)).

\begin{table}[!htb]
	\centering
	\scriptsize{	\caption{The Knapsack problem with $W=3000$}
		\label{table-2}
		\begin{tabular}{|c|cc|ccc|ccc|ccc|}
			\cline{2-12}
			\multicolumn{1}{c|}{} & \multicolumn{2}{c|}{$\mathcal U_{\text{conv}}(\mathcal S)$}  & \multicolumn{3}{c|}{$\mathcal U_\text{box}$}  & \multicolumn{3}{c|}{$\mathcal U_{\textsc{pca}}(\mathcal S, m)$}&  \multicolumn{3}{c|}{$\mathcal U_\cap (\mathcal S, m)$}   \\ \hline
			\multirow{2}[0]{*}{$\rho$} & \multirow{2}[0]{*}{Value} & Time   & \multirow{2}[0]{*}{Value} & Time & Gap     & \multirow{2}[0]{*}{Value}& Time & Gap  &  \multirow{2}[0]{*}{Value}& Time & Gap \\
			&  & (secs)   &  & (secs) & (\%)    &  & (secs) & (\%) &  & (secs) & (\%)\\\hline
			-0.8  & 7134 & 747.1   & 6249 & 0.1 & 12.41      & 6468 & 547.5  &  9.34       & 6570 & 85.7 & 7.91  \\
			-0.5  & 7082  & 739.6  & 6210 & 0.1 & 12.31      & 6296 & 457.5  & 11.10   & 6392 & 144.1 & 9.74   \\
			-0.2  & 7111  & 742.6  & 6246 & 0.1 & 12.16      & 6222 & 83.0  &  12.50    &  6338 & 70.0 & 10.87  \\
			0    & 7040   & 742.0  & 6200& 0.1  & 11.93     & 6133 & 54.3   &  12.88       & 6256 & 37.9 & 11.14   \\
			0.2  & 7158 & 741.4  & 6311 & 0.1 & 11.83      & 6214  & 50.4 & 13.19    & 6350 & 24.4 & 11.29  \\
			0.5  & 7017 & 740.9  & 6196 & 0.1  & 11.70     & 6098   & 55.3 & 13.10  &  6235 & 27.9 & 11.14 \\
			0.8  & 7109 & 742.8  & 6251 & 0.1 & 12.07      & 6175   & 32.5   & 13.14  & 6308 & 20.2 & 11.27 \\
			\hline
	\end{tabular}}
	
\end{table}

\begin{table}[!htb]
	\centering
	\scriptsize{\caption{The Knapsack problem with $W=4000$}
		\label{table-1}
		\begin{tabular}{|c|cc|ccc|ccc|ccc|}
			\cline{2-12}
			\multicolumn{1}{c|}{}& \multicolumn{2}{c|}{$\mathcal U_{\text{conv}}(\mathcal S)$}  & \multicolumn{3}{c|}{$\mathcal U_\text{box}$}  & \multicolumn{3}{c|}{$\mathcal U_{\textsc{pca}}(\mathcal S, m)$}&  \multicolumn{3}{c|}{$\mathcal U_\cap (\mathcal S, m)$}   \\ \hline
			\multirow{2}[0]{*}{$\rho$} & \multirow{2}[0]{*}{Value} & Time   & \multirow{2}[0]{*}{Value}  & Time & Gap    & \multirow{2}[0]{*}{Value} & Time & Gap  & \multirow{2}[0]{*}{Value}& Time & Gap \\
			&  & (secs)   &  & (secs) & (\%)     & & (secs) & (\%)  & & (secs) & (\%) \\ \hline
			-0.8  & 8585  & 725.2 & 7686   & 0.1 & 10.47      & 8071   & 52.57 & 5.99        & 8147  & 12.5 & 5.10\\
			-0.5  & 8484  & 728.4  & 7599  & 0.1 & 10.43      & 7852   & 30.2 & 7.45       & 7921 & 18.2 & 6.64  \\
			-0.2  & 8443  & 743.0  & 7585 & 0.1  & 10.16      & 7692  & 5.3   & 8.89      & 7780 & 19.5 & 7.85  \\
			0    & 8443   & 730.3  & 7609 & 0.1   & 9.88      & 7648  & 3.5  & 9.42       & 7746 & 11.5 & 8.26   \\
			0.2  & 8497 & 729.9  & 7646 & 0.1  & 10.02      & 7691  & 7.2   & 9.49      & 7795 & 12.5 & 8.26  \\
			0.5  & 8264  & 728.3 & 7438  & 0.1 & 10.00      & 7483   & 5.0  & 9.45       & 7581 & 12.9 & 8.26 \\
			0.8  & 8406  & 729.8 & 7575 & 0.1   & 9.89      & 7615  & 2.7   & 9.41      & 7715 & 10.8 & 8.22 \\
			\hline
	\end{tabular}}
	
\end{table}

\begin{table}[!htb]
	\centering
	\scriptsize{\caption{The Knapsack problem with $W=5000$}
		\label{table-3}
		\begin{tabular}{|c|cc|ccc|ccc|ccc|}
			\cline{2-12}
			\multicolumn{1}{c|}{} & \multicolumn{2}{c|}{$\mathcal U_{\text{conv}}(\mathcal S)$}  & \multicolumn{3}{c|}{$\mathcal U_\text{box}$}  & \multicolumn{3}{c|}{$\mathcal U_{\textsc{pca}}(\mathcal S, m)$}&  \multicolumn{3}{c|}{$\mathcal U_\cap (\mathcal S, m)$}   \\ \hline
			\multirow{2}[0]{*}{$\rho$}& \multirow{2}[0]{*}{Value} & Time   & \multirow{2}[0]{*}{Value} & Time & Gap   & \multirow{2}[0]{*}{Value}  & Time & Gap & \multirow{2}[0]{*}{Value} & Time & Gap \\
			& & (secs)   &  & (secs) & (\%)     &  & (secs) & (\%)  & & (secs)  & (\%) \\\hline
			-0.8  & 9291 & 678.9   & 8610& 0.1 & 7.33      & 9086 & 2.2  &  2.21      & 9118 & 5.8 & 1.86  \\
			-0.5  & 9442  & 666.8  & 8756& 0.1 & 7.27      & 9083   & 2.4 & 3.80   & 9128 & 6.0 & 3.33   \\
			-0.2  & 9275  & 662.4  & 8602 & 0.1 & 7.26      & 8787   & 1.7 & 5.26     & 8847  & 5.3 & 4.61  \\
			0    & 9352   & 666.2  & 8694& 0.1 & 7.04      & 8813   & 1.0  & 5.76     & 8886  & 6.1 & 4.98  \\
			0.2  & 9088 & 664.8  & 8449 & 0.1 & 7.03      & 8564  & 1.4  & 5.77    & 8635& 6.7 & 4.98  \\
			0.5  & 9263 & 664.2  & 8588 & 0.1 & 7.29      & 8716   & 1.7  & 5.91   & 8788 & 8.0 & 5.13 \\
			0.8  & 9296 & 663.1  & 8631& 0.1 & 7.15      & 8740  & 1.3  & 5.98    & 8814 & 5.8 & 5.19 \\
			\hline
	\end{tabular}}
	
\end{table}

From Tables \ref{table-2} - \ref{table-3},  we have the following observations. First, as expected $\mathcal U_\text{box}$ results in the most tractable RO Knapsack problems while $\mathcal U_{\text{conv}}(\mathcal S)$ leads to the least tractable ones. This is directly related to the high dimensionality (i.e. more decision variables) of $\mathcal U_{\text{conv}}(\mathcal S)$ compared with $\mathcal U_\text{box}$.   Second, $\mathcal U_{\text{conv}}(\mathcal S)$ results in the least conservative RO problems while $\mathcal U_\text{box}$ leads to the most conservative ones. This is because $\mathcal U_{\text{conv}}(\mathcal S)$ defines the smallest uncertainty set while $\mathcal U_\text{box}$ defines the largest uncertainty set.  Similarly, $\mathcal U_\cap (\mathcal S, m)$ leads to less conservative RO problems in comparison with $\mathcal U_{\textsc{pca}}(\mathcal S, m)$. Third, RO problems using $\mathcal U_\cap (\mathcal S, m)$ and $\mathcal U_\text{box}$, respectively, have the smallest and largest Gap values, meaning they have the most and least similar performance compared to $\mathcal U_{\text{conv}}(\mathcal S)$ in terms of conservativeness. Fourth, RO Knapsack problems with any uncertainty sets are more tractable when $W$ is larger. Fifth, $\mathcal U_\cap (\mathcal S, m)$ leads to more tractable RO problems compared to $\mathcal U_{\textsc{pca}}(\mathcal S, m)$ when $W$ is smaller while $\mathcal U_{\textsc{pca}}(\mathcal S, m)$ outperforms $\mathcal U_\cap (\mathcal S, m)$ in this regard for larger values of $W$. Sixth, when $W$ is larger, RO problems with either $\mathcal U_\text{box}$, $\mathcal U_{\textsc{pca}}(\mathcal S, m)$, or $\mathcal U_\cap (\mathcal S, m)$ have smaller Gap values. In other words, $\mathcal U_{\text{conv}}(\mathcal S)$ has a less significant benefit over other uncertainty sets in terms of conservativeness when $W$ is larger.   Seventh, when scenarios are more negatively correlated, $\mathcal U_\cap (\mathcal S, m)$ and $\mathcal U_{\textsc{pca}}(\mathcal S, m)$ result in less conservative RO problems and, moreover, they have smaller Gap values.  

We reached similar conclusions for results presented in Tables \ref{budget-conv-pca-power} - \ref{table-5}. RO problems with $\mathcal U_{\textsc{pca}}(\mathcal S, m)$ are more tractable than those with $\mathcal U_{\text{conv}}(\mathcal S)$ and less tractable than those with $\mathcal U_\text{box}$. While, RO problems with $\mathcal U_{\textsc{pca}}(\mathcal S, m)$  are less conservative than RO problems with $\mathcal U_\text{box}$ and more conservative than those with $\mathcal U_{\text{conv}}(\mathcal S)$. Lastly, similar to RO problems with $\mathcal U_{\text{conv}}(\mathcal S)$, RO problems with $\mathcal U_{\textsc{pca}}(\mathcal S, m)$ are less conservative when scenarios are more negatively correlated.  

\begin{table}[!htb]
	\centering
	\scriptsize{\caption{The power grid problem with $\rho_1=0.5$}\label{budget-conv-pca-power}
		\begin{tabular}{|c|cc|ccc|ccc|}
			\cline{2-9}
			\multicolumn{1}{c|}{} &
			\multicolumn{2}{c|}{$\mathcal U_{\text{conv}}(\mathcal S)$}  & \multicolumn{3}{c|}{$\mathcal U_\text{box}$} & \multicolumn{3}{c|}{$\mathcal U_{\textsc{pca}}(\mathcal S, m)$}  \\ \hline
			\multirow{2}[0]{*}{$\rho_2$} & Value & Time & Value & Time &  Gap   & Value & Time & Gap    
			\\ 
			&($\times10^7$)  & (secs)  & ($\times10^7$) & (secs) & (\%)    & ($\times10^7$)& (secs) & (\%)\\\hline
			-0.8  & 0.79   & 960.7  & 1.77  & 0.3 & 55.37     & 1.37  & 73.2   & 42.34       \\
			-0.5  & 0.82   & 955.6  & 1.76& 0.3   & 53.41      & 1.38 & 939.2 & 40.58        \\
			-0.2  & 0.85   & 954.5  & 1.77 & 0.3   & 51.98      & 1.48  & 816.5  & 42.57       \\
			0  & 0.88  & 956.3  & 1.77 & 0.3  & 50.28      & 1.53  & 344.1   & 42.48      \\
			0.2  & 0.88 & 961.2  & 1.77 & 0.3  & 50.28      & 1.50   & 509.3  & 41.33    \\
			0.5  & 0.94  & 962.9  & 1.75 & 0.3   & 46.29    & 1.50  & 236.6  & 37.33      \\
			0.8  & 0.93  & 1005.1  & 1.79  & 0.4 & 48.04      & 1.54  & 57.1  & 39.61     \\
			\hline
	\end{tabular}}
	
\end{table}

Based on results from Tables \ref{table-4} and \ref{table-5}, we observed the following.  First, when $\rho_1$ is larger, $\mathcal U_{\textsc{pca}}(\mathcal S, m_1)$ results in more tractable RO problems while $\mathcal U_\text{box}$ leads to less tractable RO problems. Second, RO problems with $\mathcal U_\text{box}$ are more conservative when $\rho_1$ is larger. Third, when $\rho_2$ is smaller, RO problems with $\mathcal U_{\textsc{pca}}(\mathcal S, m_1)$ have larger Gap values. In other words, RO problems with $\mathcal U_{\textsc{pca}}(\mathcal S, m_1)$ have a more remarkable benefit over those with $\mathcal U_\text{box}$ in terms of conservativeness when $\rho_2$ is smaller.  Fourth, a smaller $m_1$ results in more tractable and less conservative (i.e. robust) RO problems compared with $\mathcal U_{\textsc{pca}}(\mathcal S, m_1)$ because it yields a more smaller uncertainty set.  Therefore, $m_1$ can be used as a tool to trade-off tractability, conservativeness, and robustness of RO models with $\mathcal U_{\textsc{pca}}(\mathcal S, m_1)$.

\begin{table}[!htb]
	\centering
	\scriptsize{\caption{The power grid problem with $\rho_1=0.5$}
		\label{table-4}
		\begin{tabular}{|c|cc|ccc|ccc|ccc|}
			\cline{2-12}
			\multicolumn{1}{c|}{} &
			\multicolumn{2}{c|}{$\mathcal U_\text{box}$}  & \multicolumn{3}{c|}{$\mathcal U_{\textsc{pca}}(\mathcal S, m_1=m=48)$} & \multicolumn{3}{c|}{$\mathcal U_{\textsc{pca}}(\mathcal S, m_1=42)$}  & \multicolumn{3}{c|}{$\mathcal U_{\textsc{pca}}(\mathcal S, m_1=36)$}  \\ \hline
			\multirow{2}[0]{*}{$\rho_2$} & Value & Time & Value& Time & Gap   & Value & Time & Gap     & Value & Time & Gap 
			\\ 
			&($\times10^7$)  & (secs)  & ($\times10^7$)& (secs) & (\%)   & ($\times10^7$)& (secs) & (\%)     & ($\times10^7$)& (secs)  & (\%)  \\\hline
			-0.8  & 1.77  & 0.3 & 1.37 & 73.2 & 22.60     & 1.34 & 15.2 &  24.29      & 1.28 & 2.9 & 27.68   \\
			-0.5  &  1.76 & 0.3 &1.38 &938.2  & 21.59    & 1.35   & 63.2 & 23.30      & 1.30  & 10.6 & 26.14 \\
			-0.2  & 1.77  & 0.3 & 1.48 &816.5   & 16.38    & 1.46    & 61.3  & 17.51     & 1.42  & 9.8 & 19.77    \\
			0    &  1.77 &  0.3 & 1.53& 344.1 & 13.56    & 1.50   & 62.0  & 15.25     & 1.46    & 9.5 & 17.51     \\
			0.2   & 1.77  & 0.3  &1.50 & 509.3 & 15.25   & 1.47  & 167.9 & 16.95       & 1.43  & 34.0  & 19.21       \\
			0.5  & 1.75  &0.3  & 1.50 &236.6  & 14.29     & 1.48   & 123.4 & 15.43       & 1.43   & 41.8   & 18.29    \\
			0.8 & 1.79   & 0.4  & 1.54&57.1   & 13.97     & 1.52  & 47.7  & 15.08      & 1.49    & 12.9  & 16.76      \\
			\hline
	\end{tabular}}

\end{table}

\begin{table}[!htb]
	\centering
	\scriptsize{\caption{The power grid problem with $\rho_1=0.9$}	
		\label{table-5}
		\begin{tabular}{|c|cc|ccc|ccc|ccc|}
			\cline{2-12}
			\multicolumn{1}{c|}{}	&
			\multicolumn{2}{c|}{$\mathcal U_\text{box}$}  & \multicolumn{3}{c|}{$\mathcal U_{\textsc{pca}}(\mathcal S, m_1=m=48)$} & \multicolumn{3}{c|}{$\mathcal U_{\textsc{pca}}(\mathcal S, m_1=42)$}  & \multicolumn{3}{c|}{$\mathcal U_{\textsc{pca}}(\mathcal S, m_1=36)$}  \\\hline 
			\multirow{2}[0]{*}{$\rho_2$}& Value  & Time & Value& Time & Gap   & Value  & Time  & Gap  & Value& Time & Gap \\
			& ($\times10^7$) & (secs) 	& ($\times10^7$)& (secs) & (\%)   & ($\times10^7$)& (secs) & (\%)     & ($\times10^7$)& (secs) & (\%) \\\hline
			-0.8 &  1.86 & 0.4 & 1.31 &7.7   & 29.57    & 1.30  & 2.1 & 30.11    & 1.29 & 1.2 & 30.65   \\
			-0.5  & 1.87    & 0.4 & 1.45 & 6.0  & 22.46   & 1.44 & 1.7 & 22.99    & 1.43 & 1.4 & 23.53  \\
			-0.2 &  1.87 & 0.4 & 1.47  & 17.1  & 21.39     & 1.45  & 3.6 & 22.46     & 1.44 & 1.5  & 22.99   \\
			0    & 1.88  & 0.4 & 1.59  & 6.1  & 15.43     & 1.58 & 2.6  & 15.96      & 1.56 & 1.4  & 17.02  \\
			0.2   & 1.91  & 0.4  & 1.62 & 6.9   & 15.18   & 1.61 & 2.9  & 15.71    & 1.60 & 1.5 & 16.23    \\
			0.5 & 1.91  &  0.4 &1.63  & 2.8    & 14.66   & 1.62  & 2.5 & 15.18     & 1.61 & 1.4 & 15.71   \\
			0.8 &  1.86 &  0.4  & 1.64  & 4.1  & 11.83     & 1.64 & 9.2 & 11.83    & 1.63  & 2.1 & 12.37  \\
			\hline
	\end{tabular}}
\end{table}	

Base on to numbers in the ``Value" columns of Tables \ref{table-7} - \ref{table-8}, the power grid problem with $\mathcal U_{\textsc{pca}}(\mathcal S, m_1)$ leads to a lower bound for the power grid problem with $\mathcal U_{\textsc{pca}}(\mathcal S, m)$.  Time and gap results show that for smaller $m_1$, the power grid problem with $\mathcal U_{\textsc{pca}}(\mathcal S, m_1)$ results in a looser lower bound and is more computationally tractable. Moreover, larger values of $\rho_1$ lead to tighter and more tractable lower bounds. Again, the results justify that $m_1$ can be used as a tool to trade-off tractability, conservativeness, and robustness of RO models with $\mathcal U_{\textsc{pca}}(\mathcal S, m_1)$.  

\begin{table}[!htb]
	\centering
	\scriptsize{\caption{Lower bounds for the power grid problem with $\rho_1=0.5$}
		\label{table-7}
		\begin{tabular}{|c|cc|ccc|ccc|}
			\cline{2-9}
			\multicolumn{1}{c|}{} &
			\multicolumn{2}{c|}{$\mathcal U_{\textsc{pca}}(\mathcal S, m_1=m)$}  & \multicolumn{3}{c|}{$\mathcal U_{\textsc{pca}}(\mathcal S, m_1=42)$}  & \multicolumn{3}{c|}{$\mathcal U_{\textsc{pca}}(\mathcal S, m_1=36)$}  \\ \hline
			\multirow{2}[0]{*}{$\rho_2$} & Value & Time   & Value & Time & Gap     & Value & Time & Gap 
			\\ 
			&($\times10^7$)  & (secs)  & ($\times10^7$)& (secs) & (\%)     & ($\times10^7$)& (secs)  & (\%)  \\\hline
			-0.8  & 1.37  & 73.2 &  1.34 & 15.2 &  2.19      & 1.28 & 2.9 & 6.57   \\
			-0.5  & 1.38  & 938.2    & 1.35   & 63.2 & 2.17      & 1.30  & 10.6 & 5.80  \\
			-0.2  & 1.48  & 816.5  & 1.46    & 61.3  & 1.35     & 1.42  & 9.8 & 4.05   \\
			0    & 1.53  & 344.1 & 1.50   & 62.0  & 1.96     & 1.46    & 9.5 & 4.58    \\
			0.2   & 1.50  & 509.3  & 1.47  & 167.9 & 2.00       & 1.43  & 34.0  & 4.67      \\
			0.5  & 1.50  & 236.6 & 1.48   & 123.4 & 1.33        & 1.43   & 41.8   & 4.67    \\
			0.8 & 1.54  & 57.1   & 1.52  & 47.7  & 1.30      & 1.49    & 12.9  & 3.25      \\
			\hline
	\end{tabular}}

\end{table}

\begin{table}[!htb]
	\centering
	\scriptsize{\caption{Lower bounds for power grid problem with $\rho_1=0.9$}	
		\label{table-8}
		\begin{tabular}{|c|cc|ccc|ccc|}
			\cline{2-9}
			\multicolumn{1}{c|}{}	&
			\multicolumn{2}{c|}{$\mathcal U_{\textsc{pca}}(\mathcal S, m_1=m)$}   & \multicolumn{3}{c|}{$\mathcal U_{\textsc{pca}}(\mathcal S, m_1=42)$}  & \multicolumn{3}{c|}{$\mathcal U_{\textsc{pca}}(\mathcal S, m_1=36)$}  \\\hline 
			\multirow{2}[0]{*}{$\rho_2$}& Value  & Time    & Value  & Time  & Gap  & Value& Time & Gap \\
			& ($\times10^7$) & (secs) 	 & ($\times10^7$)& (secs) & (\%)     & ($\times10^7$)& (secs) & (\%) \\\hline
			-0.8 & 1.31  & 7.7     & 1.30  & 2.1 & 0.76    & 1.29 & 1.2 & 1.53    \\
			-0.5  & 1.45   & 6.0    & 1.44 & 1.7 & 0.69    & 1.43 & 1.4 & 1.38   \\
			-0.2 & 1.47  & 17.1     & 1.45  & 3.6 & 1.36     & 1.44 & 1.5  & 2.04   \\
			0    & 1.59  & 6.1 & 1.58 & 2.6  & 0.63      & 1.56 & 1.4  & 1.89  \\
			0.2   & 1.62  & 6.9    & 1.61 & 2.9  & 0.62    & 1.60 & 1.5 & 1.23    \\
			0.5 & 1.63  & 2.8   & 1.62  & 2.5 & 0.61     & 1.61 & 1.4 & 1.23   \\
			0.8 & 1.64  & 14.1 & 1.64 & 9.2 & 0.00    & 1.63  & 2.1 & 0.61  \\
			\hline
	\end{tabular}}
	
\end{table}

\section{Conclusion}\label{Conclusion}
In this paper, we proposed a systematic approach to develop data-driven polyhedral uncertainty sets using PCA. These uncertainty sets alleviate some of the drawbacks of conventional polyhedral uncertainty sets. Primarily, the proposed uncertainty sets capture the correlation information between uncertain parameters and are less conservative than the box uncertainty sets. Moreover, they lead to more computationally tractable RO models compared to the convex hull of uncertainty data.  The number of the leading principal components in these uncertainty sets can be used as a tool to trade-off tractability, conservativeness, and robustness of RO models. Additionally, we developed a theoretical bound on the gap between the optimal value of a static RO problem under a piece-wise linear objective function with a scenario-induced uncertainty set and that of its lower bound to quantify the quality of the lower bound. We also derived probabilistic guarantees for the performance of the proposed uncertainty sets by developing explicit lower bounds on the number of scenarios required to construct uncertainty sets.

We can extend the current research by addressing some its current limitations. First, it would be worthwhile to leverage other machine learning techniques to improve the proposed scenario-induced uncertainty sets.  Second, it would be useful to derive a theoretical bound on the gaps between the optimal value of an ARO problem, under a more general objective function, with a scenario-induced uncertainty set and that of its lower bound. Third, we used the first $m_1$ principal directions with the largest variance to develop approximate scenario-induced uncertainty sets. However, these principal directions may not always lead to the best results. Future studies may focus on developing a systematic approach to finding the $m_1$ directions that result in the best performance, with respect to computational tractability and robustness of the solution. Finally, we would like to improve the lower bound on the number of scenario samples required to achieve the desired probabilistic performance guarantees.

\bibliographystyle{apalike} 
\SingleSpacedXI
\setlength\bibsep{7pt}
\bibliography{MyReferences} 




\end{document}